\documentclass[10pt]{article}
\usepackage[a4paper, margin=2cm]{geometry}
\usepackage[fleqn]{amsmath}
\usepackage{amssymb, amsfonts, mathtools}
\usepackage{booktabs}
\usepackage{listings}
\usepackage{nccmath}
\usepackage{authblk}

\usepackage{breakcites}
\usepackage[breaklinks=true]{hyperref}  
\usepackage{doi}

\usepackage{caption}
\usepackage[parfill]{parskip}

\usepackage{sectsty}
\sectionfont{\fontsize{11}{13}\selectfont}
\subsectionfont{\fontsize{11}{13}\selectfont}

\usepackage{amsthm}
\usepackage{xargs}                      
\usepackage[pdftex,dvipsnames]{xcolor} 

\usepackage[section]{algorithm}
\usepackage{algorithmic}

\usepackage[capitalise]{cleveref}
\crefformat{equation}{(#2#1#3)} 

\usepackage{graphicx}
\graphicspath{{images/}}
\usepackage[colorinlistoftodos,prependcaption,textsize=small,color=Apricot]{todonotes}

\usepackage[utf8]{inputenc}
\usepackage[T1]{fontenc} 

\usepackage{accents}
\renewcommand{\tilde}{\widetilde}
\renewcommand{\hat}{\widehat}

\newcommand{\mc}{\mathcal}

\newcommand{\proj}{\mathcal{P}\hspace{-1.5pt}}

\DeclareMathOperator{\R}{\mathbb{R}}
\DeclareMathOperator{\E}{\mathbb{E}}

\DeclareMathOperator{\tensor}{\otimes}
\DeclareMathOperator{\rank}{rank}

\newtheorem{theorem}{Theorem}[section]
\newtheorem{lemma}[theorem]{Lemma}
\newtheorem{remark}[theorem]{Remark}
\crefname{remark}{Remark}{Remarks}

\newtheorem{proposition}[theorem]{Proposition}
\newtheorem{definition}[theorem]{Definition}

\title{Streaming Tensor Train Approximation}

\begin{document}

\author[1]{Daniel Kressner}
\author[2]{Bart Vandereycken}
\author[2]{Rik Voorhaar}
\affil[1]{\footnotesize Institute of Mathematics, EPFL, CH-1015 Lausanne, Switzerland} 
\affil[2]{\footnotesize Section of Mathematics, University of Geneva, CH-1205 Geneva, Switzerland} 
\date{}

\maketitle

\begin{abstract}
 \noindent Tensor trains are a versatile tool to compress and work with high-dimensional data and functions. In this work we introduce the Streaming Tensor Train Approximation (STTA), a new class of algorithms for approximating a given tensor $\mathcal T$ in the tensor train format. STTA accesses $\mathcal T$ exclusively via two-sided random sketches of the original data, making it streamable and easy to implement in parallel -- unlike existing deterministic and randomized tensor train approximations. This property also allows STTA to conveniently leverage structure in $\mathcal T$, such as sparsity and various low-rank tensor formats, as well as linear combinations thereof. When Gaussian random matrices are used for sketching, STTA is admissible to an analysis that builds and extends upon existing results on the generalized Nystr\"om approximation for matrices. Our results show that STTA can be expected to attain a nearly optimal approximation error if the sizes of the sketches are suitably chosen. A range of numerical experiments illustrates the performance of STTA compared to existing deterministic and randomized approaches.
\end{abstract}

\section{Introduction}

This work proposes and analyzes a new randomized algorithm for compressing a tensor $\mc T \in \R^{n_1 \times n_2 \times \cdots  \times n_d}$ of order $d$ in the so-called tensor train (TT) format~\cite{oseledetsTensorTrainDecomposition2011a}. The TT format (also called matrix product state -- MPS -- in physics \cite{schollwockDensitymatrixRenormalizationGroup2011}) represents a tensor through contractions of third-order tensors, the so-called TT cores. Offering a number of key advantages, the TT format belongs to the most popular low-rank tensor formats. In particular, its multilinear structure and its close connection to low-rank matrix factorizations allow one to leverage tools from (numerical) linear algebra, such as the singular value decomposition (SVD), in the design and analysis of algorithms. During the last decade, the TT format has shown its utility in a wide variety of applications in scientific computing and data analysis; see \cite{grasedyckLiteratureSurveyLowrank2013, bachmayrTensorNetworksHierarchical2016,  khoromskijTensorNumericalMethods2018, hackbuschTensorSpacesNumerical2019, uschmajewGeometricMethodsLowRank2020a} for an overview.

The randomized SVD, called the Halko--Martinson--Tropp (HMT) method~\cite{halkoFindingStructureRandomness2011} in the following, is a wildly popular method for obtaining a low-rank approximation to a matrix $A\in \R^{m\times n}$. It proceeds by first generating a (random) \emph{dimension reduction matrix} (DRM) $X \in \R^{n\times r}$ with $r\ll n$. Then the product $AX$ is computed and its columns are orthonormalized by a QR factorization, which we denote by  $Q=\operatorname{orth}(AX) \in \R^{n\times r}$. A low-rank approximation of $A$ is obtained by setting $A \approx \hat A_{\mathsf{HMT}} = Q(Q^\top A)$. While very simple and effective, the HMT method is not streamable in the sense of \cite{clarksonNumericalLinearAlgebra2009} since forming $Q^\top A$ requires two passes over the input data $A$. Although the sketch $AX$ is linear in $A$, which makes it cheap to update when $A$ undergoes a change of the form $A\to A + B$, this is no longer true for $Q$ and $Q^\top A$. The generalized Nyström (GN) method, originally also devised in \cite{clarksonNumericalLinearAlgebra2009}, mitigates this issue by performing a \emph{two-sided} sketch. First, it generates two DRMs $X\in \R^{n\times r}$ and $Y\in \R^{m\times (r+\ell)}$ for some $\ell \ge 1$ such that $r+\ell \ll m$. After computing the sketches $AX$, $Y^\top A$, $Y^\top AX$, one obtains a low-rank approximation by setting
\begin{equation}\label{eq:def-gnintro}
 A \approx \hat A_{\mathsf{GN}} := AX(Y^\top AX)^\dagger Y^\top A.
\end{equation}
Here $(\cdot)^\dagger$ denotes the Moore--Penrose pseudoinverse by which multiplying is mathematically equivalent to solving a linear least-squares problem. 
As all three sketches are linear in $A$ and they dominate the cost, the GN method is streamable. Moreover, the computation of the sketches can be carried out such that it requires only one pass over the data.

In this work, we present a streaming TT approximation (STTA) -- an extension of the GN method to TTs. Our algorithm is streamable and requires only one pass over the data. Moreover, unlike most existing TT approximation methods, STTA parallelizes naturally with respect to the tensor modes, an attractive feature for high-order tensors. Note that the STTA method itself does not make any assumptions on the type of DRMs. For the case of Gaussian random DRMs, we will derive an upper bound on the expected error that generalizes the matrix results from \cite{troppPracticalSketchingAlgorithms2017, nakatsukasaFastStableRandomized2020}. In particular, we prove that the expected error is quasi-optimal and close to the one obtained via TT-SVD using only moderate oversampling.

Several different randomized approaches for TT approximation have been proposed in the literature. A straightforward approach is to replace every truncated SVD inside the TT-SVD algorithm, the standard method for TT approximation \cite{oseledetsTensorTrainDecomposition2011a}, by the HMT method. This algorithm, called TT-HMT in the following, is analyzed for Gaussian random DRMs in~\cite{huberRandomizedTensorTrain2017,wolfLowRankTensor2019} and significant savings for a (very) sparse tensor $\mc T$ are reported. In~\cite{cheRandomizedAlgorithmsApproximations2019}, the TT-HMT method is presented as well, together with an adaptive variant. In~\cite{algerTensorTrainConstruction2020}, another approach is presented that is well suited for compressing high-order derivative tensors. This variant  accesses the tensor through evaluations of the corresponding multilinear map and uses random Khatri--Rao product DRMs, that is,
(columnwise) Khatri--Rao products of Gaussian random matrices instead of unstructured Gaussian random matrices.

Neither the TT-SVD algorithm nor any of the variants mentioned above is streamable or parallel in the modes; nestedness relations among the low-rank matrix factorizations constituting the TT format represent a major obstacle for doing so. The parallel streaming TT sketching (PSTT) algorithms described in~\cite{shiParallelAlgorithmsComputing2021} overcome this obstacle by first applying the HMT method to all relevant unfoldings of $\mc T$ in parallel and restoring nestedness a posteriori, without losing any of the favorable approximation properties. The work by~\cite{shiParallelAlgorithmsComputing2021}
also discusses the use of Khatri--Rao product DRMs and proposes variants that require only one pass over the input data $\mc T$. The PSTT algorithm and its  variants inherit from the HMT method the need to compute (large) QR decompositions, which cannot be cheaply updated after a linear update of the data as explained above.

The work by~\cite{daasRandomizedAlgorithmsRounding2021} focuses on TT rounding, that is, the subsequent compression of a tensor that is already in TT format.  The standard algorithm for TT rounding is essentially the TT-SVD and thus requires to first successively orthogonalize the TT cores of the input TT format, a costly step that prevents from attaining any significant savings when using the HMT method. An algorithm called ``Randomize-then-Orthogonalize'' in~\cite{daasRandomizedAlgorithmsRounding2021}
circumvents this orthogonalization step by applying the HMT method to unfoldings of the full tensor and leveraging the TT format through the use of a TT DRM, that is, a TT with Gaussian random TT cores; see also~\cite{masolomonik}. A two-sided variant that uses the GN method with two TT DRMs instead of the HMT method is also presented. None of the algorithms presented in~\cite{daasRandomizedAlgorithmsRounding2021} are easy to parallelize with respect to the modes. Also, the formulation of the algorithms does not seem to allow for streaming and requires more than one pass over the input tensor.

The STTA method proposed in this work can be viewed as a streamable, two-sided variant of the PSST algorithm from~\cite{shiParallelAlgorithmsComputing2021}. Besides two-sided sketches of the input tensor, no other large-scale computation needs to be carried out in STTA. In particular, the QR decompositions of the PSST algorithm are avoided, which makes it cheaper to carry out updates, avoids the storage of large intermediate matrices, and allows one to easily exploit sparsity and other structure in the input tensor.
As we will explain in Section \ref{sec:TT with TT}, a particular combination of the STTA method with two TT DRMs becomes \emph{mathematically} equivalent to the two-sided variant~\cite{daasRandomizedAlgorithmsRounding2021}.

Am important advantage of the one-pass nature of STTA is that it can be efficiently used to compress (or round) a structured tensor that is given in some data sparse format. For example, we explain in Section \ref{sec:structured-sketches} in some detail the implementation for compressing Tucker, CP, and sparse tensors. Furthermore, thanks to its streamability, it is straightforward to apply STTA to a sum of all these different structures. We are not aware of an existing randomized method for which all these features have been highlighted.

\subsection{Notation and preliminaries}

We first recall important properties of the tensor train (TT) format that are needed throughout the paper. For more details and proofs, we refer to \cite{bachmayrTensorNetworksHierarchical2016, hackbuschTensorSpacesNumerical2019}.

Given a matrix $A$ or a tensor $\mathcal{T}$, its entries will be indexed as $[ i_1,\ldots, i_d]$. Furthermore, when an index $i_\mu$ is denoted as $:$, it represents all the entries for the $\mu$th mode and will result in a slice of the matrix or tensor. 

A tensor $\mathcal{T}$ of size $n_1 \times n_2 \times \cdots \times n_d$ is in the TT format if it can be written element-wise as
\begin{equation}\label{eq:TT_def_element}
 \mathcal{T}[i_1,\ldots, i_d] = \sum_{\ell_1=1}^{r_1} \cdots \!\!\sum_{\ell_{d-1}=1}^{r_{d-1}}  C_1[1,i_1,\ell_1] \, C_2[\ell_1,i_2,\ell_2] \, \cdots \, C_d[\ell_{d-1},i_d,1].
\end{equation}
The third-order tensors $C_\mu$ of size $r_{\mu-1} \times n_\mu \times r_\mu$ are the TT cores (where $r_0 = r_d = 1)$. Using the matrix $C_\mu[:, i_\mu,:]$ of size $r_{\mu-1} \times r_\mu$, the relation~\eqref{eq:TT_def_element} can be written compactly as a product of $d$ matrices (where the first and last matrices collapse to a row and column vector, respectively) as follows:
\[
 \mathcal{T}[i_1,\ldots, i_d] = C_1[1,i_1,:] \, C_2[:, i_2,: ] \,\cdots\, C_d[:,i_d,1].
\]
The TT format defines a specific way to contract the tensors $C_\mu$ which can be depicted graphically using a tensor diagram. In particular,  the contractions~\eqref{eq:TT_def_element} are represented by the diagram\\
\begin{ceqn}
\begin{equation}\label{eq:diagram_def_TT}
\begin{matrix}
\includegraphics[]{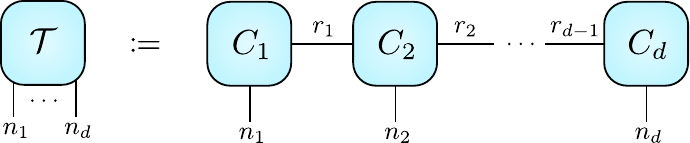}
\end{matrix}
\end{equation}
\end{ceqn}

The tuple $(r_1, \ldots, r_{d-1})$ is called the TT representation rank of the TT defined in~\eqref{eq:TT_def_element} and it determines the complexity of working with a TT. For instance, storing a tensor in  TT format requires to store the $O(dnr^2)$ entries of its TT cores, where $n:=\max_\mu(n_\mu)$ and $r:=\max_\mu(r_\mu)$. Any tensor can trivially be written in the TT format by choosing the TT representation ranks sufficiently large.

The TT representation rank of a particular tensor $\mathcal{T}$ is by no means unique but there exists an (entry-wise) minimal value which is called the TT rank of $\mathcal{T}$. The minimal value for $r_\mu$ equals the matrix rank of the $\mu$th  unfolding $\mathcal{T}$, denoted as $\mathcal{T}^{\leq \mu}$. The unfolding $\mathcal{T}^{\leq \mu}$ is one of the many ways to matricize a tensor; it is a matrix of size $(n_1 \cdots n_\mu) \times (n_{\mu+1} \cdots n_d)$ obtained from merging the first $\mu$ modes of $\mathcal T$ into row indices and the last $d-\mu$ modes into column indices.

In the rest of the paper will not distinguish between TT rank and TT representation rank and simply call $(r_1, \ldots, r_{d-1})$  the TT rank of the tensor $\mathcal{T}$ once the relation~\eqref{eq:TT_def_element} is satisfied for some cores $C_\mu$.

When splitting the tensor diagram of a TT with cores $C_\mu$ into two by cutting the edge between modes $\mu$ and $\mu+1$, one obtains two subdiagrams. Each of the two subdiagrams describes a partial contraction of the cores representing a tensor of size $n_1 \times \cdots \times n_\mu \times r_\mu$ or $n_{\mu+1} \times \cdots \times n_d \times r_\mu$. We denote the $\mu$th and $(d-\mu)$th  unfoldings of these tensors by
\[
 C_{\leq \mu} \in \mathbb{R}^{(n_1 \cdots n_\mu) \times r_\mu} \quad \text{and} \quad 
 C_{> \mu} \in \mathbb{R}^{(n_{\mu+1} \cdots n_d) \times r_\mu},
\]
which are sometimes called interface matrices.
The TT relation~\eqref{eq:TT_def_element} implies the (low-rank) factorization  $\mathcal{T}^{\leq \mu} = C_{\leq \mu}   C_{> \mu}^\top$. It also implies nestedness relations among the interface matrices
$C_{\leq \mu}$, $C_{> \mu}$ for different $\mu$:
\begin{equation} \label{eq:interfacenested}
     C_{\leq \mu} = (C_{\leq \mu-1} \otimes I) C_\mu^L  \quad \text{and} \quad 
 C_{> \mu} = (I \otimes  C_{> \mu+1} ) C_{\mu+1}^R.
\end{equation}
Here, $\otimes$ denotes the (right) Kronecker product of two matrices, and the matrices $C_\mu^L  := C_\mu^{\leq 2}$ and $C_\mu^R:= (C^{\leq 1}_\mu)^\top$ are convenient notations for the so-called left and right unfoldings of $C_\mu$, respectively.

We denote by $\| \mathcal{T} \|_F = (\sum_{i_1, \ldots, i_d} \mathcal{T}[i_1, \ldots, i_d]^2)^{1/2}$ the Frobenius norm of the tensor $\mathcal{T}$. Because of its elementwise definition, it coincides with the Frobenius norm of any matricization of $\mathcal T$; in particular $\| \mathcal{T} \|_F = \| \mathcal{T}^{\leq \mu} \|_F $ for every $\mu$.

\section{Randomized low-rank approximations} \label{sec:randlowrank}

Before describing and analyzing our algorithms for TT, we briefly review the theory of randomized methods for 
finding a low-rank approximation to a matrix $A\in \R^{m\times n}$.

As explained in the introduction, the HMT method 
constructs the approximation
\begin{equation}\label{eq:def-hmt}
  \hat A_{\mathsf{HMT}} = Q(Q^\top A),
\end{equation}
with $Q=\operatorname{orth}(AX) \in \R^{n\times r}$ for a (random) DRM $X \in \R^{n\times r}$, $r<n$. 
We let $\proj_{AX} = AX(AX)^\dagger$ denote the orthogonal projection onto the range of $AX$. Because of $\proj_{AX} = QQ^\top$, we can equivalently write \cref{eq:def-hmt} as
\begin{equation}\label{eq:hmt-proj-form}
  \hat A_{\mathsf{HMT}} = \proj_{AX} A.
\end{equation}
This approximation usually compares well with a \emph{best} rank $\hat r$ approximation $A_{\hat r}$ of slightly smaller rank $\hat r<r-1$. More specifically,
when taking expectation with respect to 
a random Gaussian\footnote{A random matrix is called Gaussian if its entries are iid standard normal random variables.} matrix $X$ one has
  \begin{equation} \label{eq:hmt}
    \E \|\hat A_{\mathsf{HMT}}-A\|_F^2 
    \le \left(1+\frac{\hat r}{r-\hat r-1}\right) \| A_{\hat r} - A \|_F^2;
  \end{equation}
  see Theorem 10.5 in~\cite{halkoFindingStructureRandomness2011}. 
  
The generalized Nyström (GN) method uses two DRMs $X\in \R^{n\times r}$ and $Y\in \R^{m\times (r+\ell)}$, for some $\ell \ge 1$ such that $r+\ell < m$, to construct the approximation
\begin{equation}\label{eq:def-gn}
  \hat A_{\mathsf{GN}} = AX(Y^\top AX)^\dagger Y^\top A.
\end{equation}
Again, this can be expressed in terms of a projection. Given two matrices $B$ and $C$ with the \textit{same number of rows}, we define the \textit{oblique projector} $\proj_{B,C} = B(C^\top B)^\dagger C^\top$. It directly follows from~\eqref{eq:def-gn} that
\begin{equation}\label{eq:gn-proj-form}
  \hat A_{\mathsf{GN}} = \proj_{AX,Y}A.
\end{equation}
The fact that $\proj_{AX,Y}$ is no longer orthogonal implies that 
the GN approximation is  (slightly) worse than the HMT approximation for the same $X$. In terms of the expected error the following bound holds. 
\begin{theorem}[{\cite[Theorem 4.3]{troppPracticalSketchingAlgorithms2017}}]\label{thm:GN}
  Consider random Gaussian matrices $X$ and $Y$ of size  $n \times r$ and $n \times (r+\ell)$, respectively, with $\ell > 0$. Then the expected approximation error of \cref{eq:def-gn} is bounded by
  \begin{equation}\label{eq:gn-error-bound}
    \E \|\hat A-A\|_F^2 
    \leq \left(1+\frac{r}{\ell-1}\right)\left(1+\frac{\hat r}{r-\hat r-1}\right)\| A_{\hat r} - A \|_F^2,
  \end{equation}
  where $A_{\hat r}$ is any best rank $\hat r$ approximation of $A$ with $\hat r<r-1$.
\end{theorem}

\begin{remark}\label[remark]{remark:gn-transpose}
The formulation of~\cref{thm:GN} assumes that $Y$ has more columns than $X$. Since
\begin{equation}
    \|(I-\proj_{AX,Y})A\|_F = \|[(I-\proj_{AX,Y})A]^\top\|_F = \|(I-\proj_{A^\top Y, X})A^\top\|_F,
\end{equation}
the bound \cref{eq:gn-error-bound} also holds in the reverse situation, that is, when $X$ has size $n\times(r+\ell)$ and $Y$ has size $n\times r$.
\end{remark}

Let us emphasize again that the GN method is \emph{streamable} and \emph{one-pass}. It consists of two steps:
\begin{enumerate}
  \item Sketch $AX$, $Y^\top A$, $Y^\top AX$.
  \item Assemble $\hat A = AX(Y^\top AX)^\dagger Y^\top A$.
\end{enumerate}
The first step is linear in $A$, and the sketches are cheap to store or communicate if $r\ll \min(m,n)$. Moreover, this step typically dominates the overall computational effort, requiring $O(rmn)$ operations for dense matrices, whereas the second step only requires $O((r+\ell)r^2)$ operations; see~\cite{nakatsukasaFastStableRandomized2020} for implementation aspects.

\section{Streaming TT approximation (STTA)}\label{sec:STTA}

We now explain our streamable TT approximation algorithm for a given tensor $\mc T$ of size $n_1\times \cdots\times n_d$.  The idea is to sketch every unfolding ${\mc T}^{\leq \mu} \in \R^{(n_1 \cdots n_\mu) \times (n_{\mu+1}\cdots n_d)}$ using several \textit{dimension reduction matrices} (DRMs) $X_\mu$, $Y_\mu$ of matching sizes:
\begin{equation}\label{eq:def_Xmu_Ymu}
\begin{aligned}
  &X_\mu\text{ is a `right' DRM of size }(n_{\mu+1}\cdots n_d) \times r_\mu^R\text{ for }1\leq \mu\leq d-1; \\
  &Y_\mu\text{ is a `left' DRM of size }(n_1\cdots n_\mu)\times r_\mu^L\text{ for }1\leq \mu\leq d-1.
\end{aligned}
\end{equation}
We require that either $r_\mu^R < r_\mu^L-1$  for all $\mu$ or $r_\mu^L < r_\mu^R-1$ for all $\mu$. Our approximation will have TT rank $r_\mu \coloneqq \min(r_\mu^L,r_\mu^R)$ for $\mu = 1,\ldots,d-1$ and is obtained from the following sketches:
\begin{equation}\label{eq:def_Psimu_Omegamu}
\begin{aligned}
  &\Psi_\mu^L = (Y_{\mu-1}^\top\tensor I) \mc T^{\leq \mu} X_\mu\text{ of size }(r_{\mu-1}^L   n_\mu) \times r_\mu^R\text{ for }1\leq \mu\leq d;\\
  &\text{$\Omega_\mu = Y_{\mu}^\top \mc T^{\leq \mu} X_\mu$ of size $r_{\mu}^L \times r_\mu^R$ for $1\leq \mu\leq d-1$.}
\end{aligned}
\end{equation}
For simplicity we use the convention that $Y_0=X_d=1$, so that \cref{eq:def_Psimu_Omegamu} remains valid for $\mu=1$ and $\mu=d$.
Recall that $\Psi_\mu^L$ denotes the left unfolding of the $r_{\mu-1}\times n_\mu \times r_\mu^R$ tensor $\Psi_\mu$. In terms of tensor diagrams, the relations~\eqref{eq:def_Psimu_Omegamu} take the following form:
\begin{ceqn}
\begin{equation}\label{eq:diagram_def_Omega_Psi}
\begin{matrix}
\includegraphics[]{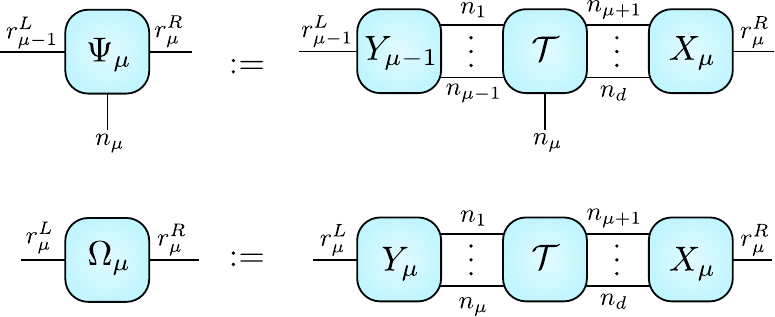}
\end{matrix}
\end{equation}
\end{ceqn}
Using these sketches, we construct an approximation $\tilde{\mc T}$ as follows:
\begin{ceqn}
\begin{equation}\label{eq:diagram_def_T}
\begin{matrix}
\includegraphics[]{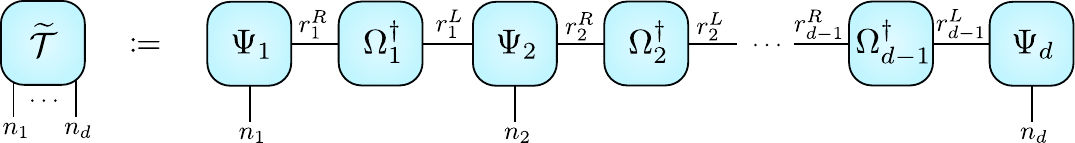}
\end{matrix}
\end{equation}
\end{ceqn}
To arrive at a tensor diagram that takes the form of a TT as in \eqref{eq:diagram_def_TT}, one needs to contract each $\Omega_\mu^\dagger$ with an adjacent order-three tensor. When $r_{\mu}^R < r_{\mu}^L$, it is preferable to choose the tensor to the right, because otherwise the (representation) ranks are unnecessarily high. More specifically, we obtain a TT representation for $\tilde{\mc T}$ by defining each
$r_{\mu-1}^R\times n_\mu \times r_\mu^R$ core $C_\mu$ via its right unfolding 
\begin{equation}\label{eq:def_Cmu_R}
C_\mu^{R} := \begin{cases}
      \Psi_1, & \mu = 1,\\
      \Omega_{\mu-1}^\dagger\Psi_\mu^R, & \mu > 1.
    \end{cases} 
\end{equation}
Assuming that $\Omega_\mu$ has full (column) rank, we can equivalently obtain  $C_\mu^{R}$ as the solution of the linear least-squares problem
\begin{equation}\label{eq:Cmu-as-lstsq}
     \min_X \|  \Omega_{\mu-1}X - \Psi_\mu^R \|_F^2.
\end{equation}
Especially when $\Omega_\mu$ is not well-conditioned, care needs to be taken when solving this problem numerically; see also~\cite[\S 4.1]{nakatsukasaFastStableRandomized2020}. 
In our implementation, we use the LAPACK routine \texttt{gelsd} which utilizes a (truncated) SVD of $\Omega_{\mu-1}$ and discards singular values smaller than $\epsilon_{\text{mach}} \| \Omega_{\mu-1}\|_2$, where $\epsilon_{\text{mach}}$ denotes machine precision.

When $r_{\mu}^L > r_{\mu}^L$ it is preferable to contract $\Omega_\mu^\dagger$ with the tensor on the left, that is, the TT core $C_\mu$ is an $r_{\mu-1}^L\times n_\mu \times r_\mu^L$ tensor defined via
\begin{equation}\label{eq:def_Cmu_L}
    C_\mu^{L} := \begin{cases}
      \Psi_\mu^L\Omega_\mu^\dagger, & \mu < d\\
      \Psi_d, & \mu = d.
    \end{cases}
\end{equation}
If $\Omega_\mu$ has full row rank, this is equivalent to solving the linear least-squares problem $\min_X \| X \Omega_\mu - \Psi_\mu^L \|_F^2$.

Like the matrix case, our streaming tensor sketch algorithm (STTA) can be split into two subroutines. First we \emph{sketch} the tensor $\mc T$, and then we \emph{assemble} the TT from the sketches. These two routines are summarized below in \cref{alg:sketch,alg:assemble}. For convenience we only state the version for $r_\mu^R < r_\mu^L$, which produces a TT of rank $(r_1^R,\dots,r_{d-1}^R)$; the other version is analogous.

\begin{algorithm}
\caption{STTA: \texttt{sketch}}\label{alg:sketch}
\begin{algorithmic}[1]
    \REQUIRE Tensor $\mathcal T$ of shape $n_1\times \dots \times n_d$; Right DRMs $X_\mu$ of shape $(n_{\mu+1} \cdots n_d) \times r_\mu^R$ and left DRMs $Y_\mu$ of shape $(n_1\cdots n_\mu) \times r_\mu^L$ for $\mu=1$ to $d-1$  with $r_\mu^R < r_\mu^L$. 
    \ENSURE Sketches $\Psi_\mu$ of shape $r_{\mu-1}^L\times n_\mu\times r_{\mu}^R$ for $1\leq \mu\leq d$; Sketches $\Omega_\mu$ of shape $r_\mu^L \times r_\mu^R$ for $1\leq \mu \leq d-1$.
    
    \FOR{$\mu=1$ to $d$}
    \STATE $\mathcal T^{\leq \mu} \leftarrow \mathtt{reshape}(\mathcal T,\, (n_1\cdots n_\mu) \times (n_{\mu+1}\cdots n_d))$
    \IF{$\mu=1$}
    \STATE $\Psi_1 \leftarrow \mathcal T^{\leq 1} X_1$
    \ELSIF{$\mu<d$}
    \STATE $\Psi_\mu \leftarrow (Y_{\mu-1}^\top\tensor I_{n_\mu}) \mathcal T^{\leq \mu} X_\mu$
    \ELSE
    \STATE $\Psi_d \leftarrow (Y_{d-1}^\top\tensor I_{n_d}) \mathcal T^{\leq d} $
    \ENDIF
    \IF{$\mu<d$}
    \STATE $\Omega_\mu \leftarrow Y_\mu^\top \mathcal T^{\leq \mu} X_\mu$
    \ENDIF
    \ENDFOR
\end{algorithmic}
\end{algorithm}

\begin{algorithm}
\caption{STTA: \texttt{assemble}}\label{alg:assemble}
\begin{algorithmic}[1]
    \REQUIRE Sketches $\Psi_\mu$ of shape $r_{\mu-1}^L\times n_\mu\times r_{\mu}^R$ for $1\leq \mu\leq d$; Sketches $\Omega_\mu$ of shape $r_\mu^L \times r_\mu^R$ for $1\leq \mu \leq d-1$.
    \ENSURE Cores $C_\mu$ of shape $r_{\mu-1}^R\times n_\mu\times r_\mu^R$ for $1\leq \mu\leq d$ of a tensor train approximation of $\mc T$.
    
    \STATE $C_1 \leftarrow \mathtt{reshape}(\Psi_1,1\times n_1\times r_1^R)$
    \FOR{$\mu=2$ to $d-1$}
        \STATE $B_\mu^R \leftarrow \mathtt{lstsq}(\Omega_{\mu-1},\,\Psi_\mu^R)$\hfill // Stable least-square solver for $\Omega_{\mu-1}^\dagger \Psi^R_\mu$ as in \cref{eq:Cmu-as-lstsq}
        \STATE $C_\mu \leftarrow \mathtt{reshape}(B_\mu^R,\,r_{\mu-1}^R\times n_\mu\times r_\mu^R)$
    \ENDFOR
\end{algorithmic}
\end{algorithm}

\paragraph{DRMs:} Algorithm~\ref{alg:assemble} takes $2d-2$ right and left DRMs as input. The simplest choice for DRMs are Gaussian random matrices, which lead to matrices $X_\mu$ and $Y_\mu$ that are independent with respect to  $\mu$. When $\mathcal{T}$ is itself in low-rank format (TT, CP, and Tucker formats -- see Sections~\ref{sec:cp} and~\ref{sec:tucker}), it is computationally advantageous to use a single TT DRM (a TT with Gaussian
random TT cores) to define all right DRMs, and similarly for all left DRMs. In turn, the right DRMs $X_1, \ldots X_{d-1}$ and the left DRMs $Y_1, \ldots Y_{d-1}$ are no longer independent for different $\mu$ but coupled through the nestedness relations~\eqref{eq:interfacenested}; see Definition~\ref{def:tt-drm} for details.

For the special case of TT DRMs, the approximation returned by STTA is mathematically equivalent to the one produced by Algorithm~3.3 of~\cite{daasRandomizedAlgorithmsRounding2021}; see also Section \ref{sec:TT with TT}.

\paragraph{Streaming:} Because \cref{alg:sketch} is linear in $\mc T$, it is indeed  \textit{streamable} and \textit{one-pass}. Suppose that $\mc T$ is decomposed as a sum $\mc
T_1+\cdots +\mc T_k$ and we compute the sketches $\Psi_{\mu,i}$, $\Omega_{\mu,i}$ of each
summand $\mc T_i$. The sketch of the sum is then the sum of the sketches: $\Psi_\mu = \sum_{i=1}^k \Psi_{\mu,i}$, and similarly
for $\Omega_\mu$. These sketches can be performed in a data-parallel fashion, and only
one round of communication is required, involving a total $O(nr^2dk)$ of data\footnote{This is assuming that the DRMs are not communicated. In practice, the DRMs are always pseudo-random so that they can be efficiently generated locally on each node using a seeded method. This will be discussed further in \cref{sec:structured-sketches}},
making it very suitable for a distributed setting.

\paragraph{Parallelism:} In addition to the data-parallelism due to streaming, we note that all for loops in both \cref{alg:sketch,alg:assemble} can be executed in parallel. This is also true for TT DRMs at the expense of some redundant computations. Moreover, the most expensive operations in \cref{alg:sketch} are  matrix multiplications, which enjoy a high degree of parallelism themselves.

\paragraph{Adaptive rank and block sketching:} We can adaptively increase the rank of the approximation by performing a new sketch. Suppose we have computed a rank-$r$ approximation of $\mc T$ using the sketches $\Psi_\mu^1$, $\Omega_\mu^1$, and DRMs $X_\mu^1$, $Y_\mu^1$ with $1\leq\mu\leq d$. We can then compute a rank-$(r+r')$ approximation by first generating a second pair of DRMs $X_\mu^2$, $Y_\mu^2$ of rank $r'$. Sketching $\mc T$ with the rank-$(r+r')$ DRMs $X_\mu = [X_\mu^1\ X_\mu^2]$ and $Y_\mu = [Y_\mu^1\ Y_\mu^2]$ is mathematically equivalent to performing the following block-wise sketches:
\begin{equation}\label{eq:block-sketch}
  \begin{split}
    \Psi_\mu &= \begin{bmatrix}
      ((Y_{\mu-1}^1)^\top\tensor I) \mc T^{\leq \mu} X_\mu^1 & ((Y_{\mu-1}^1)^\top\tensor I) \mc T^{\leq \mu} X_\mu^2 \\
      ((Y_{\mu-1}^2)^\top\tensor I) \mc T^{\leq \mu} X_\mu^1 & ((Y_{\mu-1}^2)^\top\tensor I) \mc T^{\leq \mu} X_\mu^2
    \end{bmatrix} \\ 
    \Omega_\mu &= \begin{bmatrix}
      (Y_{\mu}^1)^\top \mc T^{\leq \mu} X_\mu^1 & (Y_{\mu}^1)^\top \mc T^{\leq \mu} X_\mu^2 \\
      (Y_{\mu}^2)^\top \mc T^{\leq \mu} X_\mu^1 & (Y_{\mu}^2)^\top \mc T^{\leq \mu} X_\mu^2
    \end{bmatrix}.
  \end{split}
\end{equation}

Note that in this scenario the top-left blocks of both sketches are $\Psi_\mu^1$ and $\Omega_\mu^1$, respectively, which were both already computed. This means that first computing a rank-$r$ sketch, and then computing a rank-$(r+r')$ sketch is not more expensive than directly computing a rank-$(r+r')$ sketch. Furthermore, \cref{eq:block-sketch} gives a block algorithm for computing the sketch, which further increases the available parallelism.

\subsection{Error analysis}

For the error analysis of STTA, it is beneficial to express the obtained approximation in terms of projectors, like for the HMT method in~\cref{eq:hmt-proj-form} and for the GN method in~\cref{eq:gn-proj-form}. For this purpose, let us define the oblique projectors
\begin{equation}\label{eq:def_Pmu}
    \proj_\mu := \proj_{\mc T^{\leq \mu} X_\mu, Y_\mu} = \mc T^{\leq \mu} X_\mu (Y_\mu^\top\mc T^{\leq \mu} X_\mu)^\dagger Y_\mu^\top = \mc T^{\leq \mu} X_\mu\Omega_\mu^\dagger Y_\mu^\top.
\end{equation}
We claim that the approximation \cref{eq:diagram_def_T} of STTA satisfies 
\begin{equation}\label{eq:approximation-proj-form}
    \tilde {\mc T}^{\leq d-1}=(\proj_1\tensor I)\cdots (\proj_{d-2}\tensor I)\proj_{d-1}\mc T^{\leq d-1},
\end{equation}
where the size of the identity matrices $I$ is such that all matrix products are well defined. In terms of a tensor diagram, this expression takes the form
\begin{ceqn}
\begin{equation}\label{eq:diagram_def_T_with_proj}
\begin{matrix}
\includegraphics{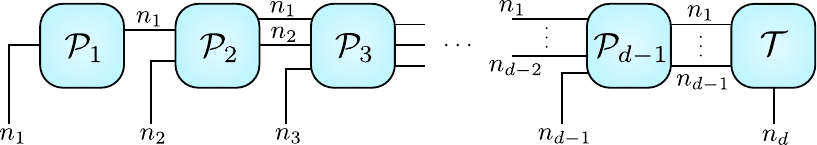}
\end{matrix}
\end{equation}
\end{ceqn}
To show~\eqref{eq:approximation-proj-form}, we first note that  the definitions \cref{eq:def_Psimu_Omegamu} and \cref{eq:def_Pmu} imply
\[
\left(\left[\Psi_1^L\Omega_1^\dagger\right] \tensor I\right)\Psi_2^L 
= \left(\left[ \mc T^{\leq 1} X_1 \left(Y_1^\top \mc T^{\leq 1}X_1\right)^\dagger Y_1^\top \right] \tensor I\right) \mc T^{\leq 2}X_2 = (\proj_1 \tensor I)\mc T^{\leq 2}X_2,
\]
which corresponds to the diagram
\begin{ceqn}
  \begin{equation*}
  \begin{matrix}
  \includegraphics[]{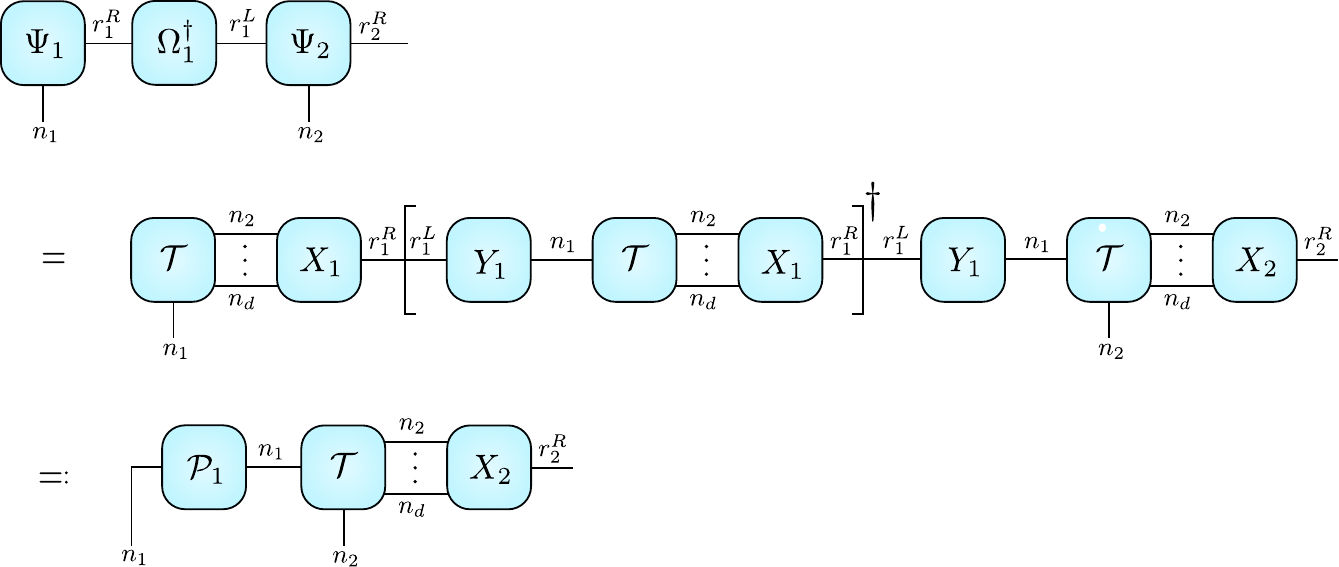}
  \end{matrix}
  \end{equation*}
\end{ceqn}
Comparing with~\eqref{eq:diagram_def_T}, this shows that the left part of the network representing $\tilde{\mc T}$ can be replaced by   $\mc P_1$. More generally, we have the relation 
\[
\left(\left[\mc T^{\leq \mu}X_\mu\Omega_\mu^\dagger\right] \tensor I\right)\Psi_{\mu+1}^L 
= \left(\left[ \mc T^{\leq \mu} X_\mu \left(Y_\mu^\top \mc T^{\leq \mu}X_\mu\right)^\dagger Y_\mu^\top \right] \tensor I\right) \mc T^{\leq \mu+1}X_{\mu+1} = (\proj_\mu \tensor I)\mc T^{\leq \mu+1}X_{\mu+1},
\]
which allows us to successively replace, from the left to the right, the network representing $\tilde {\mc T}$ with the projectors $\proj_\mu$, 
until we arrive at
\begin{ceqn}
  \begin{equation*}
  \begin{matrix}
  \includegraphics[]{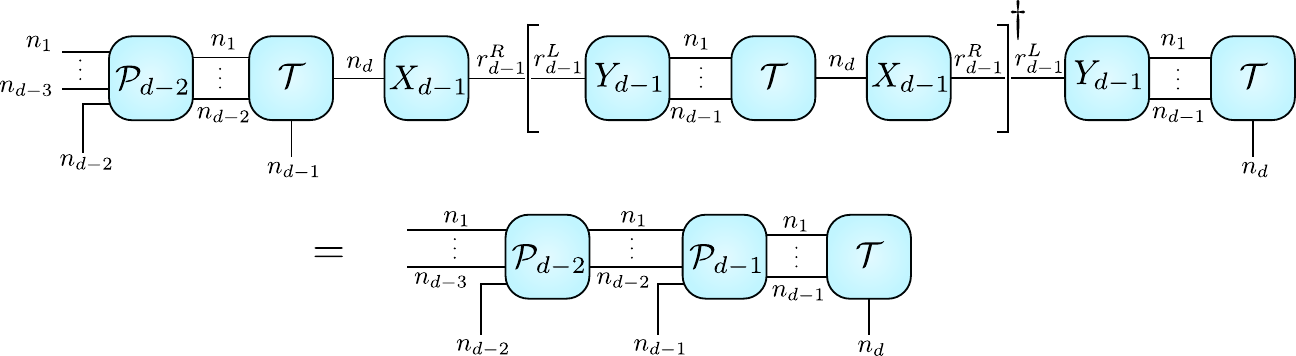}
  \end{matrix}
  \end{equation*}
\end{ceqn}
In turn, we obtain the network~\cref{eq:diagram_def_T_with_proj} as claimed.

\begin{proposition}\label[proposition]{prop:proj-estimate}
  The approximation $\tilde{\mc T}$ returned by STTA satisfies\footnote{Note that the non-commuting product $\textstyle{\prod}_{\alpha < \mu} (\proj_\alpha\tensor I)$ in \cref{eq:prop-proj-estimate} is evaluated in natural increasing order with respect to $\alpha$. Furthermore, for $\mu=1$ the product is empty and corresponds to the identity matrix.}
\begin{equation}\label{eq:prop-proj-estimate}
    \|\tilde{\mc T}-\mc T\|_F \leq \sum_{\mu=1}^{d-1}\left\|\textstyle{\prod}_{\alpha < \mu} (\proj_\alpha\tensor I)(I-\proj_{\mu}){\mc T}^{\leq \mu}\right\|_F.
  \end{equation}
\end{proposition}
\begin{proof}
  Using~\cref{eq:approximation-proj-form}, we have that 
  \begin{align*}
   \|\mc T - \tilde{\mc T}\|_F &=  \|[I-(\proj_1\otimes I)\cdots (\proj_{d-2}\otimes I)\proj_{d-1})]{\mc T}^{\leq d-1}\|_F \\ 
    &\leq \|[I-(\proj_1\otimes I)\cdots (\proj_{d-2}\otimes I)]{\mc T}^{\leq d-1}\|_F+\|(\proj_1\otimes I)\cdots (\proj_{d-2}\otimes I)[I-\proj_{d-1}]{\mc T}^{\leq d-1}\|_F. 
  \end{align*}
We can now continue the argument inductively since
\[
 \|[I-(\proj_1\otimes I)\cdots (\proj_{d-2}\otimes I)]{\mc T}^{\leq d-1}\|_F = \|[I-(\proj_1\otimes I)\cdots \proj_{d-2}]{\mc T}^{\leq d-2}\|_F. \qedhere
 \]
\end{proof}

\paragraph{Gaussian DRMs.}

We now focus on the particular case when Gaussian DRMs $X_\mu, Y_\mu$ are used and extend the bounds from Section~\ref{sec:randlowrank} to the STTA algorithm.

\begin{theorem}\label{thm:error-bound}
  Suppose that the DRMs $X_\mu,Y_\mu$ defined in \cref{eq:def_Xmu_Ymu} are independent standard Gaussian matrices with $r_\mu^R < r_\mu^L-1$ and that the rank of $\mathcal{T}^{\le \mu}$ is not smaller than $r_\mu^R$ for $\mu = 1,\ldots,d-1$. 
Then for any $(\hat r_1,\dots,\hat r_{d-1})$ such that $\hat r_\mu < r_\mu^R-1$, the STTA approximation $\tilde {\mc T}$ of $\mathcal T$ obtained from \cref{alg:assemble}  satisfies
  \begin{equation}\label{eq:main-error-bound}
    \E\|\tilde {\mc T}-\mc T\|_F \leq  \sum_{\mu=1}^{d-1}{\left[\textstyle\prod_{\alpha=1}^{\mu-1} c_\alpha\right]c'_{\mu} \sqrt{\sum_{k>\hat r_\mu}\sigma_k(\mc T^{\leq \mu})^2}} \leq \left( \sum_{\mu=1}^{d-1} {\left[\textstyle\prod_{\alpha=1}^{\mu-1} c_\alpha\right]c'_{\mu} } \right) \ \|\mathcal T_{\hat {\mathbf r}}-\mathcal T\|_F,
  \end{equation}
  where $\mathcal T_{\hat {\mathbf r}}$ is any best TT approximation of $\mathcal T$ of TT rank $(\hat r_1,\dots,\hat r_d)$ , and $c_\mu$, $c_\mu'$ are given by 
  \begin{equation}\label{eq:def_c_mu}
    c_{\mu} \coloneqq 1+\sqrt{ \frac{r^R_\mu}{r_\mu^L-r_\mu^R-1}},\qquad c'_{\mu} \coloneqq \sqrt{1+\frac{r_\mu^R}{r_\mu^L-r_\mu^R-1}}\cdot \sqrt{ 1+\frac{\hat r_\mu}{r_\mu^R-\hat r_\mu -1}}.
  \end{equation}
\end{theorem}
\begin{proof}
By the assumptions, both $\mc T^{\leq \mu} X_\mu$ and $Y_\mu^\top \mc T^{\leq \mu} X_\mu$ have full column rank almost surely. We consider arbitrary $1\le \mu < d-1$ and temporarily consider $X_\mu$ fixed (such that the full-rank condition holds).
The QR decomposition $\mc T^{\leq \mu} X_\mu=QR$ yields an invertible  factor $R$ of size $r_\mu^R \times r_\mu^R$, which allows us to rewrite the definition~\cref{eq:def_Pmu} of $\proj_\mu$ as follows:
\[
\proj_\mu = \mc T^{\leq \mu} X_\mu(Y_\mu^\top \mc T^{\leq \mu} X_\mu)^\dagger Y_\mu^\top = Q R(Y_\mu^\top QR)^\dagger Y_\mu^\top 
 = Q(Y_\mu^\top Q)^\dagger Y_\mu^\top.
\]
We now complete the $r_\mu^R$ columns of $Q$ to a square orthogonal matrix $\begin{bmatrix}Q & Q_\perp \end{bmatrix}$. Using $QQ^\top + Q_\perp Q_\perp^\top = I$ we obtain for any matrix $A$ of size $(n_{1}\cdots n_\mu) \times (n_{\mu+1}\cdots n_d)$ that
\begin{align*}
    \|\proj_\mu A\|_F & \le \|(Y_\mu^\top Q)^\dagger Y_\mu^\top Q Q^\top A\|_F + \|(Y_\mu^\top Q)^\dagger Y_\mu^\top Q_\perp Q_\perp^\top A\|_F \\
&= \| Q^\top A\|_F + \|(Y_\mu^\top Q)^\dagger Y_\mu^\top Q_\perp Q_\perp^\top A\|_F \le \|A\|_F + \|(Y_\mu^\top Q)^\dagger Y_\mu^\top Q_\perp Q_\perp^\top A\|_F.
\end{align*}
By the orthogonal invariance of Gaussian random matrices, $Z_1\coloneqq Y_\mu^\top Q$ 
and $Z_2 \coloneqq Y_\mu^\top Q_\perp$ are independent Gaussian random matrices of size $r_\mu^L\times r_\mu^R$ and $r_\mu^L\times (n_1\cdots n_\mu-r_\mu^R)$, respectively. This allows us to apply Propositions 10.1 and 10.2 from \cite{halkoFindingStructureRandomness2011} to obtain
\[
\E_{Z_1,Z_2} \|Z_1^\dagger Z_2 Q^\top_\perp\!A\|_F \le \sqrt{ \E_{Z_1,Z_2} \|Z_1^\dagger Z_2 Q^\top_\perp\!A\|_F^2} =
\sqrt{ \E_{Z_1}\|Z_1^\dagger\|_F^2}\,\|Q^\top_\perp\!A\|_F = \sqrt{\frac{r_\mu^R}{r_\mu^L-r_\mu^R-1}} \|Q^\top_\perp\!A\|_F.
\]
In summary, we have
\[
\E_{X_1,Y_1} \|\proj_\mu A\|_F \le c_\mu \|A\|_F.
\]

Now, by the law of total expectation, it holds that
\begin{align*}
\E_{X_\mu,Y_\mu}&{}_{,\dots,,X_\mu,Y_\mu}\|\textstyle{\prod}_{\alpha = 1}^{\mu-1} (\proj_\alpha\tensor I)(I-\proj_{\mu}) {\mc T}^{\leq \mu}\|_F \\ 
&=\E_{X_2,Y_2,\dots,X_\mu,Y_\mu}\left[\E_{X_1,Y_1}\|(\proj_1\tensor I) \textstyle{\prod}_{\alpha = 2}^{\mu-1} (\proj_\alpha\tensor I) (I-\proj_{\mu}){\mc T}^{\leq \mu}\|_F\ \big|\ X_2,Y_2,\dots,Y_{\mu},X_{\mu} \right]\\
& \le c_{1}\E_{X_2,Y_2,\dots,X_\mu,Y_\mu}\left[\|\textstyle{\prod}_{\alpha = 2}^{\mu-1} (\proj_\alpha\tensor I)   (I-\proj_{\mu}) {\mc T}^{\leq \mu}\|_F \right],
\end{align*}
where we used that only $\proj_1$ depends on $X_1,Y_1$. 
Continuing inductively we obtain
\begin{equation}
  \begin{split}
  \E_{X_1,Y_1,\dots,X_{\mu},Y_\mu}\|\textstyle{\prod}_{\alpha = 1}^{\mu-1} (\proj_\alpha\tensor I)(I-\proj_{\mu}){\mc T}^{\leq \mu}\|_F
  \le  c_{1}\cdots c_{\mu-1} \E_{X_\mu,Y_\mu}\|(I-\proj_\mu)\mc T^{\leq \mu}\|_F.
  \end{split}
\end{equation}
The last factor coincides with the GN approximation error, which -- according to  \cref{thm:GN} -- is bounded by
\begin{align*}
    \E_{X_\mu,Y_\mu}\|(I-\proj_\mu) {\mc T}^{\leq \mu}\|_F^2 &\leq \left(1+\frac{r_\mu^R}{r_\mu^L-r_\mu^R-1}\right)\left(1+\frac{\hat r_\mu}{r_\mu^R-\hat r_\mu -1}\right)\sum_{k>\hat r_\mu} \sigma_k (\mc T^{\leq \mu})^2\\
    &= (c'_{\mu})^2 \sum_{k>\hat r_\mu} \sigma_k (\mc T^{\leq \mu})^2
\end{align*}
for any $\hat r_{\mu}<r_\mu^R-1$, where $\sigma_k(\cdot)$ denotes the $k$th larges singular value of a matrix.
 Finally, we apply \cref{prop:proj-estimate} to obtain the first inequality of \cref{eq:main-error-bound}. The second inequality follows directly from the general bound
\begin{equation}\label{eq:sing_val_lower_bound_best_rank}
  \sum_{k>\hat r_\mu} \sigma_k (\mc T^{\leq \mu})^2 \leq \|\mathcal T-\mathcal T_{\hat{\mathbf r}}\|_F^2;
\end{equation}
see, e.g., \cite[Thm.~11.6]{hackbuschTensorSpacesNumerical2019}. \end{proof}

Observe that if one uses $r_\mu^L=\kappa r_\mu^R$ for a constant $\kappa>1$ in \cref{thm:error-bound}, then $c_\mu=O(1)$. As $\kappa$ increases, $c_\mu$ approaches one, which mitigates the exponential dependence on $d$ in the constants.

The condition $r_\mu^R < r_\mu^L-1$ for all $\mu$ of \cref{thm:error-bound} can be replaced by $r_\mu^L < r_\mu^R-1$, using a modification of the proof that proceeds in the other direction, from `right to left'. The approximation bound~\cref{eq:main-error-bound} still holds but with the roles of $r_\mu^R$ and $r_\mu^L$ reversed in~\cref{eq:def_c_mu}.

\paragraph{Random TT DRMs.} The proof of 
\cref{thm:error-bound} heavily relies  on the DRMs being Gaussian. As mentioned before and demonstrated in~\cite{daasRandomizedAlgorithmsRounding2021, masolomonik}, 
it can often be computationally beneficial to use TT DRMs, for which we now give a formal definition.

\begin{definition}[TT DRM]\label[definition]{def:tt-drm}
Consider a tensor in TT format with random TT cores $C_\mu \in \R^{r^L_{\mu-1}\times n_\mu \times r^L_\mu}$ for which the entries are i.i.d. normal random variables with zero mean and variance $1/r^L_{\mu}$. Then \emph{left TT DRMs}  $(Y_1,\ldots,Y_{d-1})$, required as input by Algorithm~\ref{alg:sketch}, are defined via the (left) interface matrices
\[
Y_\mu := C_{\le \mu} \in \R^{(n_1\cdots n_\mu) \times r_\mu}, \quad \mu = 1,\ldots, d-1, 
\]
which satisfy the nestedness relations~\eqref{eq:interfacenested}, that is, $Y_\mu = (Y_{\mu-1} \otimes I) C^L_\mu$.

Analogously, \emph{right TT DRMs}
$(X_1,\ldots,X_{d-1})$ are defined through  $X_\mu:=C_{> \mu}$ for another (independent) set of random TT cores $C_\mu \in \R^{r^R_{\mu-1}\times n_\mu \times r^R_\mu}$ with variance $1/r^R_{\mu-1}$.
\end{definition}

The normalization of the variance in Definition~\ref{def:tt-drm} is chosen so that we have $\E \| X C_\mu^L\|_F^2 = \|X\|_F^2$ for each TT core $C_\mu$ and any fixed matrix $X$ of appropriate size. This avoids numerical difficulties, such as over- or underflow, in particular when dealing with tensors of high order as in \Cref{sec:experiment-tensor-order}. 

Proving an analogue to \cref{thm:error-bound} for TT DRMs is difficult. As a first, preliminary result, we prove that the product $\proj = (\proj_1\tensor I)\cdots (\proj_{d-2}\tensor
I)\proj_{d-1}$, which produces the STTA approximation $\tilde{\mathcal T}=\proj \, \mc T$, is itself an oblique projector. This intriguing property does not hold for Gaussian DRMs.

\begin{lemma}\label[lemma]{lem:nested-ranges}
  Suppose that $r^R_\mu \ge r^L_\mu$ for all $\mu$, and $(Y_1,\ldots,Y_{d-1})$ is a left TT DRM. If $\rank(\mc T^{\leq \mu}) \geq r_\mu^R$ for all $\mu$, 
  then $\proj:= (\proj_1\tensor I)\cdots (\proj_{d-2}\tensor I)\proj_{d-1}$ is an oblique projector almost surely.
\end{lemma}

\proof We first establish $Y_\mu(\proj_{\mu-1}\tensor I)=Y_\mu$ from the nestedness relations~\eqref{eq:interfacenested}:
\begin{align*}
  Y^\top_\mu(\proj_{\mu-1}\tensor I) &= (C^L_\mu)^\top(Y_{\mu-1}^\top\tensor I)[\mc T^{\leq \mu-1} X_{\mu-1}(Y_{\mu-1}^\top \mc T^{\leq \mu-1} X_{\mu-1})^\dagger Y_{\mu-1}^\top\tensor I]\\ 
  &=  (C^L_\mu)^\top[Y_{\mu-1}^\top\mc T^{\leq \mu-1} X_{\mu-1}(Y_{\mu-1}^\top \mc T^{\leq \mu-1} X_{\mu-1})^\dagger Y_{\mu-1}^\top\tensor I\ \\
  &= (C^L_\mu)^\top (Y_{\mu-1}^\top \tensor I) = Y_\mu.
\end{align*}
Here, we used that the $r_\mu^L\times r_\mu^R$ matrix $Y_{\mu-1}^\top \mc T^{\leq \mu-1} X_{\mu-1}$ has full row rank almost surely because of $r_\mu^R \ge r_\mu^L$. It now follows that
\begin{equation}\label{eq:prod-two-proj}
  \proj_\mu(\proj_{\mu-1}\tensor I) =  \mc T^{\leq \mu} X_\mu\Omega_\mu^\dagger Y_\mu^\top (\proj_{\mu-1}\tensor I) =  \mc T^{\leq \mu} X_\mu\Omega_\mu^\dagger Y_\mu^\top = \proj_\mu.
\end{equation}

We now show that $(\proj_1\tensor I)\cdots (\proj_{d-2}\tensor I)\proj_{d-1}$ is a projector. For convenience, we omit writing the $\tensor I$ factors. Observe that by applying \cref{eq:prod-two-proj} multiple times, we get the identities
\[
  \proj_\mu \cdots \proj_{2} \proj_{1} = \proj_{\mu} \cdots\proj_{3} \proj_{2} =  \dots = \proj_\mu
\]
and
\[
  (\proj_\mu\cdots \proj_2\proj_{1})(\proj_1\proj_{2}\cdots \proj_{\mu}) = \proj_\mu\cdots \proj_2\proj_1\proj_{2}\cdots \proj_{\mu} = 
  \proj_\mu\cdots \proj_3\proj_2\proj_{3}\cdots \proj_{\mu} = \cdots =
  \proj_\mu,
\]
which allow us to show that
\[
  (\proj_1\cdots \proj_\mu)(\proj_1\cdots\proj_\mu) = (\proj_1\cdots \proj_\mu)(\proj_\mu\cdots\proj_1)(\proj_1\cdots\proj_\mu) = (\proj_1\cdots \proj_\mu),
\]
as required.\qed

The assumption $r_\mu^R \ge r_\mu^L$ can be replaced by  $r_\mu^L \ge r_\mu^R$. In this case, the STTA approximation satisfies  $\tilde{\mathcal T}=\mc T\proj$ for an oblique projector $\proj$ associated with right TT DRMs ($X_1,\ldots,X_{d-1}$).

\section{Sketching structured tensors}\label{sec:structured-sketches}

We recall that STTA requires the sketches $\Psi_\mu$, $\Omega_\mu$  defined in \cref{eq:diagram_def_Omega_Psi}. For a general order-$d$ tensor $\mc T \in \R^{n_1 \times \cdots  \times n_d}$, the straightforward computation of these sketches requires a cost of $O(d n^d)$, which becomes excessive for large $d$. On the other hand, large-scale tensors usually have some \textit{structure} that allows for faster sketching. In this section, we will explain how the STTA algorithm can be efficiently implemented for tensors with the following structures:
\begin{itemize}
  \item sparse tensors;
  \item TT tensors (recompression);
  \item CP tensors;
  \item Tucker tensors.
\end{itemize}
Among these, only sparse tensors can be efficiently sketched using Gaussian DRMs;  for the other formats we will instead  use TT DRMs for sketching. Although our error bounds only apply to Gaussian DRMs, we will see that TT DRMs deliver similar performance.

Because of the linearity of sketching, we automatically obtain efficient algorithms for approximating linear combinations of any of these types of structured tensors. For example, one might be interested in sketching a sum of a sparse tensor and a TT.

\subsection{Sketching sparse tensors}

Suppose we want to approximate an $n_1\times\cdots\times n_d$ sparse tensor $\mathcal T$ with $N \ll n_1\cdots n_d$ nonzero entries. These entries are represented via a list $\mathcal E=(\mathcal E_1,\dots,\mathcal E_N)\in\mathbb R^N$ of values and a list $\mathcal J$ of $N$ multi-indices $(i_j^{(1)},\dots,i_j^{(d)})$, with $1\leq i_j^{(\mu)} \leq n_\mu$ and $1\leq j\leq N$. Given arbitrary DRMs $X_\mu$, $Y_\mu$, the sketches $\Omega_\mu$, $\Psi_\mu$ can be computed entry-wise as follows: 
\begin{align}
  &\Omega_\mu[k,l] = (Y_\mu^\top \mathcal T^{\leq \mu} X_\mu)[k,l] = \sum_{j=1}^N \mathcal E_j Y_\mu\!\left[i_j^{(1)},\dots,i_j^{(\mu)},k\right]X_\mu\!\left[i_j^{(\mu+1)},\dots,i_j^{(d)},l\right]\label{eq:sparse-tensor-contract-omega} \\ 
  &\Psi_\mu[k,l,m] = \left((Y_{\mu-1}^\top\tensor I)\mathcal T^{\leq \mu} X_\mu\right)\![k,l,m] = 
  \sum_{\substack{j=1\\ i_{j}^{(\mu)}=l}}^N \mathcal E_j Y_{\mu-1}\!\left[i_j^{(1)},\dots,i_j^{(\mu-1)},k\right]X_\mu\!\left[i_j^{(\mu+1)},\dots,i_j^{(d)},m\right].\label{eq:sparse-tensor-contract-psi}
\end{align}
Note that we tacitly reshaped $Y_\mu$ into an $n_1\times \dots n_\mu\times r_\mu^L$ tensor in order to access its rows by the multi-index $(i_j^{(1)},\dots,i_j^{(\mu)})$, and similarly for $X_\mu$ and $Y_{\mu-1}$. Furthermore, $\Psi_\mu$ is viewed as a $r_{\mu-1}^L \times n_\mu \times r_\mu^R$ tensor.

\paragraph{Gaussian DRMs.}

For larger $d$ and modest $N$,
most rows of $Y_\mu$ and $X_\mu$  are never needed in \cref{eq:sparse-tensor-contract-omega}. This allows us to efficiently use Gaussian DRMs in this context if we only generate those rows of $X_\mu$ and $Y_\mu$ that will actually be accessed. In order to do this in a distributed setting, we need to generate rows of a DRM \textit{on demand} and \emph{consistently} for each mode $\mu$. One such way is to use a simple hashing algorithm to convert the indices $(i^{(1)}_j,\dots,i^{(\mu)}_j)$ indexing the row of the DRM into a pseudorandom integer, which is then cast into a floating point number between 0 and 1, and finally converted into a normally distributed number. In the context of sketching this seems to produce pseudorandom numbers that are as good as those generated by more advanced pseudorandom number generators. 

The formal procedure is detailed below in \cref{alg:hash-rng}. The seed $s$ is specific to each mode $\mu$ and ensures that the entries of $X_\mu$ (or $Y_\mu$) are consistent in a distributed setting. The uniform numbers are converted to normally distributed numbers by applying the inverse CDF. For the hashing function, we used  \texttt{splitmix64} which generates adequate pseudo-random numbers.\footnote{See \url{https://rosettacode.org/wiki/Pseudo-random_numbers/Splitmix64} and \url{https://docs.oracle.com/javase/8/docs/api/java/util/SplittableRandom.html}.}

\renewcommand{\algorithmiccomment}[1]{\hfill// #1}
\begin{algorithm}
  \caption{Hashed pseudorandom number generator that returns one row of a Gaussian DRM.}\label{alg:hash-rng}
  \begin{algorithmic}[1]
    \REQUIRE Indices $(i_1,\dots,i_\mu)$ and integers $(n_1,\dots,n_\mu)$ such that $1 \le i_j \leq n_j$, integer $r\geq 1$, seed $s$.
    \ENSURE $(x_1,\dots,x_r)$ pseudorandom numbers with i.i.d. normal distribution 
    \STATE $m\leftarrow 0$
    \FOR{$j=1$ to $\mu$} 
      \STATE $m\leftarrow n_j (m + i_j)$ \COMMENT{Turn multi-index into linear index}
    \ENDFOR
    \STATE $m\leftarrow mr + \mathtt{hash}(s)$
    \FOR{$j=1$ to $r$}
      \STATE $h_j \leftarrow \mathtt{hash}(m+j)$ \COMMENT{$\mathtt{hash}$ produces unsigned integers}
      \STATE $x_j \leftarrow$ Interpret $h_j$ as floating point number\footnotemark and extract mantissa.
      \STATE $x_j \leftarrow 2x_j-1$.
      \STATE $x_j\leftarrow \mathtt{normal\_inverse\_cdf}(x_j)$ \COMMENT{Apply inverse CDF of the normal distribution}
    \ENDFOR
  \end{algorithmic}
\end{algorithm}
\footnotetext{We first set the leftmost 3 bits of $h_j$ to \texttt{001} to avoid NaN/infinite flags, and ensure a nonzero exponent. This does not affect the mantissa. Recall that the mantissa of a floating point number with nonzero exponent lies in the interval $[0.5,1]$.}

Using \cref{alg:hash-rng} we only generate \textit{at most} $N$ rows of each DRM, requiring a total cost of $O(Nrd)$ to generate all DRMs. The asymptotic cost of computing $\Omega_\mu$ using \cref{eq:sparse-tensor-contract-omega} is therefore $O(Nr^2)$ flops. Computing $\Psi_\mu$ using \cref{eq:sparse-tensor-contract-psi} has the same cost. This gives a total asymptotic cost for all the sketches with Gaussian DRMs of $O(Nr^2d)$, with straightforward parallelization with respect to the number of nonzero entries $N$.

\paragraph{TT DRMs.}
TT DRMs can also be used to efficiently sketch sparse tensors. This is because of the nestedness relations~\eqref{eq:interfacenested}, which allow us to efficiently retrieve the entries with indices $[i_j^{(1)},\dots,i_j^{(\mu)},k]$ of the left interface matrices $C^{\le \mu}$ subsequently for $\mu = 1,\ldots,d$. \cref{alg:gather} describes the resulting procedure. This  allows us to compute within $O(Nr^2d)$ operations all the entries of the TT DRMs $Y_\mu$ and $X_\mu$ needed for sketching in one sweep for all $\mu$. With these entries at hand, we can sketch with TT DRMs using \cref{eq:sparse-tensor-contract-omega} and \cref{eq:sparse-tensor-contract-psi} at the same asymptotic cost as when sketching with Gaussian DRMs.

\begin{algorithm}
  \caption{Selected entries of left interface matrices for a TT}\label{alg:gather}
  \begin{algorithmic}[1]
    \REQUIRE Tensor train cores $C_1,\ldots,C_d$ with each $C_\mu$ of size $r_{\mu-1}\times n_\mu\times r_\mu$;\\
    List of indices $\mc J = \{(i_1^{(1)},\dots,i_1^{(d)}),\dots,(i_N^{(1)},\dots,i_N^{(d)})\}$.
    \ENSURE Matrices $V_1,\dots,V_{d}$ of sizes $N\times r_\mu$ with entries of left interface matrices $C^{\leq \mu}$ corresponding to $\mc J$.

    \STATE $V_1 \leftarrow \mathbf 0_{N\times r_1}$
    \FOR{$j=1$ to $N$}
      \STATE $V_1[j,\colon] \leftarrow C_1[1,\,i_j^{(1)},\colon]$
    \ENDFOR
    \FOR{$\mu=2$ to $d$}
      \STATE $V_\mu \leftarrow \mathbf 0_{N\times r_\mu}$
      \FOR{$\ell=1$ to $n_\mu$}
        \STATE $M\leftarrow C_\mu[\colon,\ell,\colon]$
        \FOR{$j=1$ to $N$ such that $i_j^{(\mu)}=\ell$}
          \STATE $V_\mu[j,\colon] \leftarrow V_{\mu-1}[j,\colon]M$
        \ENDFOR
      \ENDFOR
    \ENDFOR

  \end{algorithmic}
\end{algorithm}

\subsection{Sketching tensor trains}\label{sec:TT with TT}

Let us now consider a tensor $\mathcal T$ in TT format with a  relatively high TT rank $(s_1,\ldots,s_{d-1})$, which we aim to compress to lower TT rank $(r_1,\ldots,r_{d-1})$, that is, $r_\mu < s_\mu$. 
It is clear that Gaussian DRMs are not well suited for this purpose simply because all entries of the DRMs need to be accessed at least once during the sketch computation, yielding a complexity of at least $O(n^{d})$ for computing the sketches. On the other hand, we will see that sketching with TT DRMs yields a computational cost that is linear in $d$ (assuming fixed TT ranks).

Let $B_\mu$ denote the cores of a random TT of TT rank $(r_1^L,\dots,r_{d-1}^L)$ that defines a left TT DRM $(Y_1,\ldots, Y_{d-1})$ according to~\cref{def:tt-drm}. Similarly, let
 $A_\mu$ denote the cores of a random TT of TT rank $(r_1^R,\dots,r_{d-1}^R)$ that defines a right TT DRM $(X_1,\ldots, X_{d-1})$. Denoting the cores of $\mathcal T$ by $C_\mu$ and inserting the TT diagrams for $X_\mu$, $Y_\mu$, $\mathcal T$ into the tensor diagrams~\cref{eq:diagram_def_Omega_Psi}, we obtain the following diagrams for $\Psi_\mu$ and $\Omega_\mu$:
\begin{ceqn}
  \begin{equation*}
      \begin{matrix}
        \includegraphics[]{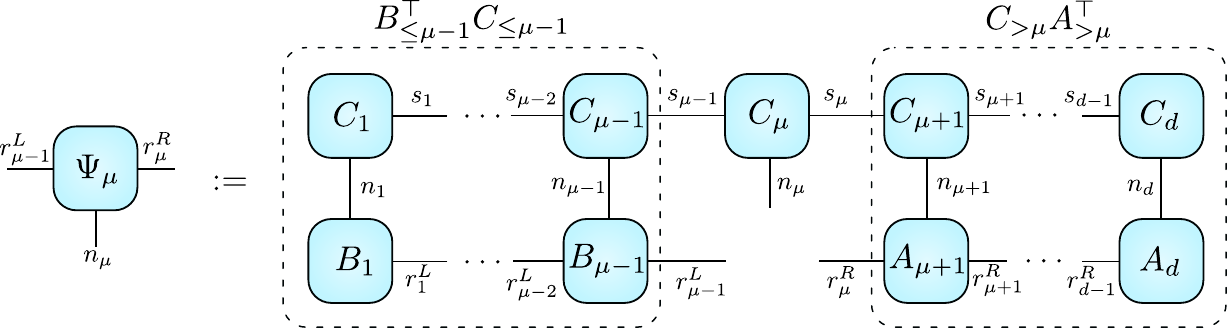} \\ \vspace{0.3cm}\\
        \includegraphics[]{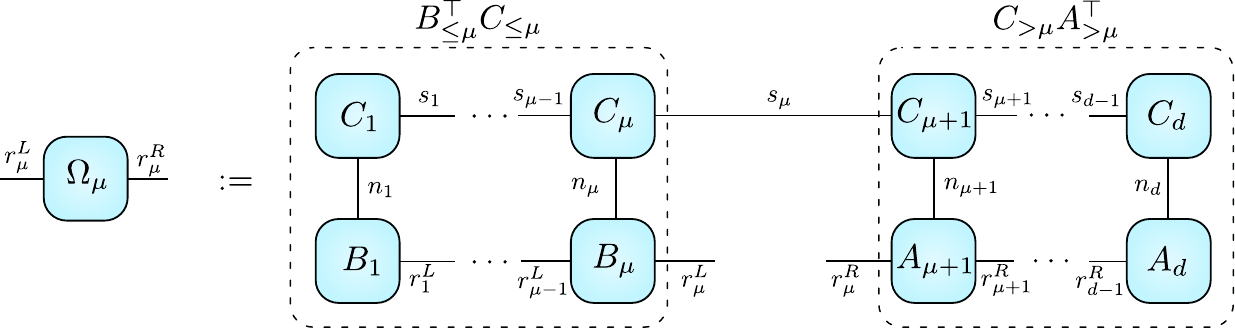}
      \end{matrix}
  \end{equation*}
\end{ceqn}

Let $r = \max_\mu(r^L_\mu, r^R_\mu)$ and $s = \max_\mu s_\mu$. All the products $B_{\leq \mu}^\top C_{\leq \mu}$ and $C_{>\mu}A_{>\mu}^\top$ (indicated by boxes above) can be computed efficiently for every $1\leq \mu\leq d-1$ in a single sweep, at a total cost of $O(rs^2 n d)$ flops. Given $B_{\leq \mu}^\top C_{\leq \mu}$ and $C_{>\mu}A_{>\mu}^\top$, we can compute each $\Psi_\mu$ at a cost of $O(rs^2 n)$ flops, and $\Omega_\mu$ at a cost of $O(r^2s)$ flops. This gives a total complexity of $O(rs^2 nd)$ for sketching a TT using TT DRMs. The procedure is summarized in \cref{alg:sketch-tt}, an efficient implementation of \cref{alg:sketch} when TT DRMs are used to sketch a TT tensor  $\mathcal T$.

\begin{algorithm}
  \caption{STTA: \texttt{sketch} a TT tensor with TT DRMs}\label{alg:sketch-tt}
  \begin{algorithmic}[1]
    \REQUIRE TT to be sketched with cores $C_1,\dots,C_d$; Left TT DRM with cores $B_1,\dots,B_d$;\\ Right TT DRM with cores $A_1,\dots,A_d$.
    \ENSURE  Sketches $\Psi_1,\dots,\Psi_d$, and $\Omega_1,\dots,\Omega_{d-1}$. 

    \STATE $L_1\leftarrow B_1^\top C_1$
    \FOR{$\mu=2$ to $d-1$}
      \STATE $L_\mu[\colon,\colon] \leftarrow \sum_{ijk} L_{\mu-1}[i,\,j]B_{\mu}[i,\,k,\colon]C_{\mu}[j,\,k,\colon]$
    \ENDFOR
    \STATE $R_{d-1}\leftarrow C_dA_d^\top$
    \FOR{$\mu=d-2$ to $2$}
      \STATE $R_{\mu}[\colon,\colon] \leftarrow \sum_{ijk} R_{\mu+1}[i,\,j]C_{\mu+1}[\colon,\,k,\,i]A_{\mu+1}[\colon,\,k,\,j]$
    \ENDFOR
    \FOR{$\mu=1$ to $d-1$}
      \STATE $\Omega_\mu \leftarrow L_\mu R_\mu$
    \ENDFOR
    \STATE $\Psi_1\leftarrow C_1R_1$
    \FOR{$\mu=2$ to $d-1$}
      \STATE $\Psi_\mu[\colon,\colon,\colon] \leftarrow \sum_{ij}L_\mu[\colon,\, i] C_\mu[i,\colon,\,j] R_\mu[j,\colon]$
    \ENDFOR
    \STATE $\Psi_d\leftarrow L_dC_d$
  \end{algorithmic}
\end{algorithm}

We remark that \cref{alg:sketch-tt} can be adapted for sketching matrix product operators (MPO)~\cite{Orus2014}, an extension of the TT format to linear operators. This is because an MPO can be transformed into a TT through reshaping its cores.

\begin{remark}
\cref{alg:sketch-tt} partly reproduces the two-sided method, Algorithm~3.3, from~\cite{daasRandomizedAlgorithmsRounding2021}. In particular, both algorithms compute the matrices
$L_\mu = B_{\le \mu}^\top C_{\le \mu}$ and 
$R_\mu = C_{> \mu} A_{> \mu}^\top$ in the same fashion. There are however also differences. In the notation of this paper, Algorithm 3.3 from~\cite{daasRandomizedAlgorithmsRounding2021} computes a factored pseudo-inverse of $\Omega_\mu = L_\mu R_\mu$ and then uses these factors to compress the cores $C_\mu$ of the input tensor. This last step has the same complexity as the overall algorithm; it is however not linear in $C_\mu$ and thus not streamable.
While mathematically equivalent, our algorithms (\cref{alg:sketch-tt} combined with \cref{alg:assemble}) avoid the explicit computation of the factored pseudo-inverse. Through the use of $\Psi_\mu$, the majority of the computation (\cref{alg:sketch-tt}) is linear in each $C_\mu$ and thus streamable. \cref{alg:assemble} also does not involve the (potentially large) cores of the input tensor.
\end{remark}

\subsection{Sketching canonical polyadic decompositions} \label{sec:cp}

The canonical polyadic (CP) decomposition~\cite{Kolda2009} represents a tensor $\mc T$ as a sum of $N$ rank-1 tensors.
The CP decomposition can be expressed in terms of factor matrices $V_1 \in \R^{N\times n_1},\dots,V_d \in \R^{N\times n_d}$ as follows:
\begin{equation}
  \mc T = \sum_{i=1}^N V_1[i,:]\tensor V_2[i,:]\tensor\cdots\tensor V_d[i,:].
\end{equation}
The tensor diagram corresponding to this equation is given by 
\begin{ceqn}
  \begin{equation*}
      \begin{matrix}
        \includegraphics[]{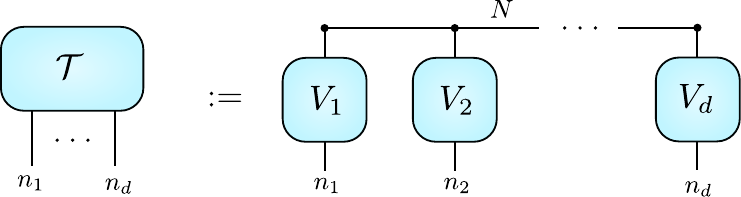}
      \end{matrix}
  \end{equation*}
\end{ceqn}

As in the TT case, sketching $\mc T$ with Gaussian DRMs would be too expensive and we therefore resort again to left/right TT DRMs with cores $A_1,\dots,A_d$ and $B_1,\dots,B_d$, respectively. The sketch $\Omega_\mu=Y_\mu^\top \mc T^{\leq \mu} X_\mu$ is obtained from
\begin{equation}\label{eq:omega-cp-sketch}
  \Omega_\mu = L_\mu R_\mu \coloneqq \big[B_{\leq \mu}^\top \left(V_1\odot \cdots \odot  V_\mu\right)^\top\big]\big[ \left(V_{\mu+1}\odot \cdots\odot V_d\right)A_{> \mu}^\top\big],
\end{equation}
where $\odot$ denotes the row-wise Khatri-Rao product, that is, 
$V_{\mu} \odot V_{\mu+1}$ is the $N\times n_\mu n_{\mu+1}$ matrix that contains in the $i$th row all products $V_{\mu}[i_1,i_2]V_{\mu+1}[i_1,i_3]$. 
The relation~\eqref{eq:omega-cp-sketch} in terms of diagrams:
\begin{ceqn}
  \begin{equation*}
      \begin{matrix}
        \includegraphics[]{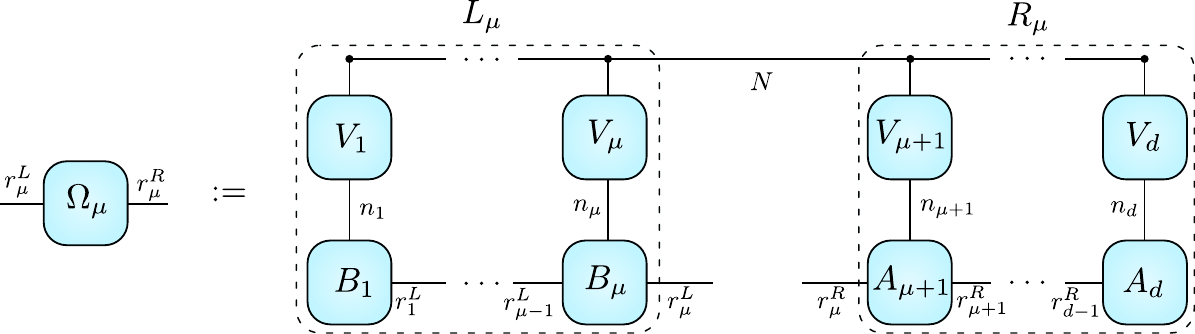}
      \end{matrix}
  \end{equation*}
\end{ceqn}

Similarly we obtain the sketch $\Psi_{\mu}$ using the following contraction
\begin{equation}\label{eq:psi-cp-sketch}
  \Psi_\mu[i_1,i_2,i_3] = \sum_j L_{\mu-1}[i_1,j]V_\mu[j,i_2]R_\mu[j,i_3],
\end{equation}
or in terms of diagrams:
\begin{ceqn}
  \begin{equation*}
      \begin{matrix}
        \includegraphics[]{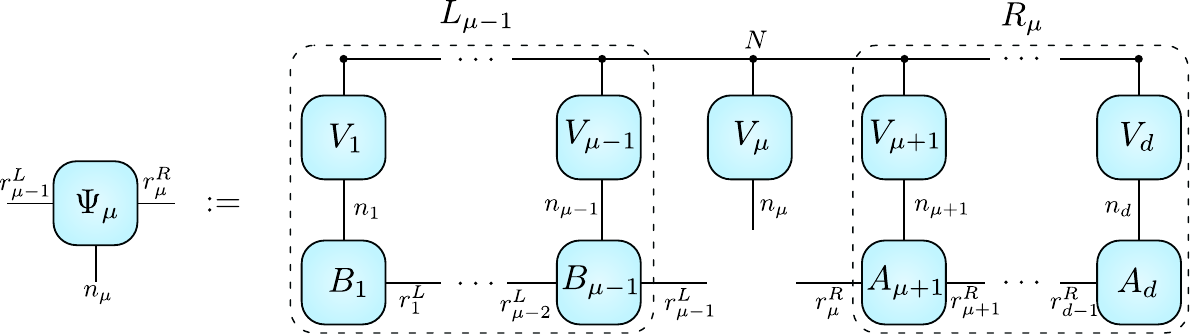}
      \end{matrix}
  \end{equation*}
\end{ceqn}

The matrices $L_\mu$ and $R_\mu$ can be efficiently computed in a single sweep for all $1\leq \mu\leq d$ at a total cost of $O(Nnr^2d)$ flops. Given $L_\mu$ and $R_\mu$, then $\Omega_\mu$ can be computed at a cost of $O(Nr^2)$ flops, and $\Psi_\mu$ at a cost of $O(Nnr^2)$ flops. This gives a total cost of $O(Nnr^2d)$ to sketch a CP tensor using TT DRMs. The procedure is summarized in \Cref{alg:sketch-cp}. 

\begin{algorithm}
  \caption{STTA: \texttt{sketch} a tensor in CP decomposition with TT DRMs}\label{alg:sketch-cp}
  \begin{algorithmic}[1]
    \REQUIRE CP to be sketched with factor matrices $V_1,\dots,V_d$, left TT DRM with cores $B_1,\dots,B_d$,\\ right TT DRM with cores $A_1,\dots,A_d$.
    \ENSURE  Sketches $\Psi_1,\dots,\Psi_d$ and $\Omega_1,\dots,\Omega_{d-1}$. 

    \STATE $L_1\leftarrow B_1^\top V_1^\top$ 
    \FOR{$\mu=2$ to $d-1$}
      \STATE $L_\mu[i_1,i_2] \leftarrow \sum_{j_1,j_2}L_{\mu-1}[j_1,i_2] V_{\mu}[i_2,j_2] B_\mu[j_1,j_2,i_1]$
    \ENDFOR
    \STATE $R_{d-1}\leftarrow V_d A_d^\top$
    \FOR{$\mu=d-2$ to $2$}
      \STATE $R_\mu[i_1,i_2] \leftarrow \sum_{j_1,j_2}R_{\mu+1}[i_1,j_1]V_{\mu+1}[i_1,j_2]B_{\mu+1}[i_2,j_2,j_1]$
    \ENDFOR
    \FOR{$\mu=1$ to $d-1$}
      \STATE $\Omega_\mu\leftarrow L_\mu R_\mu$
    \ENDFOR
    \STATE $\Psi_1[1,i_1,i_2] \leftarrow \sum_j V_1[j,i_1]R[j,i_2]$
    \FOR{$\mu=2$ to $d-2$}
      \STATE $\Psi_\mu[i_1,i_2,i_3] = \sum_j L_{\mu-1}[i_1,j]V_\mu[j,i_2]R_\mu[j,i_3]$
    \ENDFOR
    \STATE $\Psi_d[i_1,i_2,1] \leftarrow \sum_j V_\mu[j,i_1]R_\mu[j,i_2]$.    
  \end{algorithmic}
\end{algorithm}

\subsection{Sketching Tucker decomposition}
\label{sec:tucker}

The Tucker decomposition~\cite{Kolda2009} represents a tensor $\mc T$ of size $n_1 \times \cdots \times n_d$ by a core tensor $\mc C$ of (usually smaller) size $s_1 \times \cdots \times s_d$ multiplied in each mode by a factor matrix
$U_\mu\in \R^{n_\mu \times s_\mu}$:
\begin{equation}
  \mc T[i_1,\dots,i_d] = \sum_{j_1,\dots,j_d}\mc C[j_1,\dots,j_d] U_1[i_1,j_1]U_2[i_2,j_2]\cdots U_d[i_d,j_d].
\end{equation}
If the columns of every $U_\mu$ are orthonormal, an effective way to obtain  a TT approximation of $\mc T$ is to first approximate the core $\mc C$ and then contract the resulting TT with the factor matrices; see, e.g.,~\cite{shiParallelAlgorithmsComputing2021}.
While orthonormality can always be ensured by computing QR decompositions of the factor matrices, the resulting updates of the core tensor come at a cost of $O(d s^{d+1})$ and typically destroy all sparsity structure in $\mathcal C$. Sketching $\mathcal T$ directly avoids these disadvantages and comes with further benefits; for example, its linearity allows us to combine tensors with different types of structure.

We agains use left/right TT DRMs with cores $A_1,\dots,A_d$ and $B_1,\dots,B_d$, respectively, for sketching a tensor in Tucker format.  The sketches $\Omega_\mu$ and $\Psi_\mu$ for a Tucker tensor $\mathcal T$ are computed according to the following tensor diagrams:
\[
\includegraphics[]{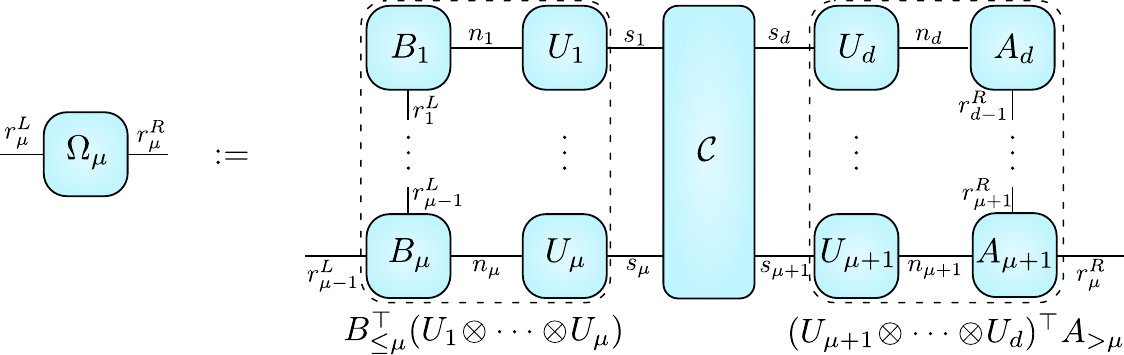}
\]

\[
\includegraphics[]{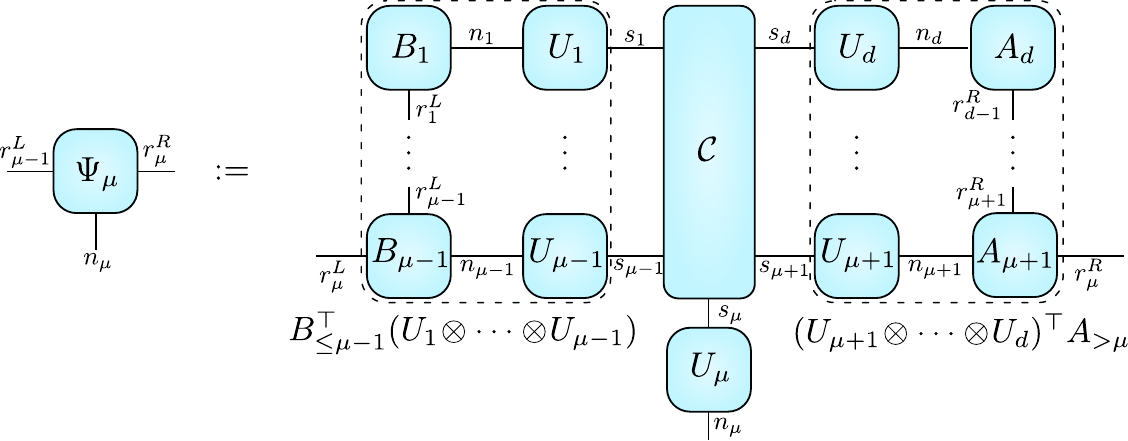}
\]

Carrying out the multiplications $B^\top_{\leq \mu}(U_1\tensor\cdots \tensor U_\mu)$ and $(U_{\mu+1}\otimes \cdots U_d)^\top A^\top_{>\mu}$  within the TT format takes $O(dn r^2s)$ flops. The subsequent explicit computation of these matrices by carrying out the contractions in a suitable order requires $O(r^2 s^{d-1})$. The cost of computing each of the sketches $\Omega_\mu$ and $\Psi_\mu$ is bounded by $O(rs^d+r^2 s^{d-1} + n r^2 s )$.
This gives a total cost for sketching the Tucker representation using TT DRMs of $O(d (n r^2s + rs^d+r^2 s^{d-1}) )$.

\section{Numerical experiments}\label{sec:numerical-experiments}

In the following, we compare STTA numerically with the classical TT-SVD algorithm and the TT-HMT method mentioned in the introduction. The TT-SVD method is used as baseline since it usually has the smallest approximation error for a fixed rank. The method is however very costly since it requires truncated SVDs of large matricizations. It is also not streamable and less suited for large structured tensors.

STTA and TT-HMT on the other hand are designed for larger tensors. While both methods have similar flop counts for \emph{dense} tensors, TT-HMT does not support streaming and is difficult to implement efficiently in a distributed setting. In its standard formulations of~\cite{huberRandomizedTensorTrain2017,wolfLowRankTensor2019, cheRandomizedAlgorithmsApproximations2019}, TT-HMT is also  not well suited for structured tensors (but we present a reformulation in the next section that partly alleviates this).

When varying the TT rank in the experiments, we increase it uniformly among all modes whenever possible. At the left/right borders of the TT, the TT ranks may be constrained by the size of the tensor: Given a prescribed rank $r$ we set $r_1 = \min\{r,n_1\}$, $r_{d-1} =  \min\{r,n_d\}$, and so on. Unless noted otherwise we set $r_\mu^L=2r_\mu^R$ for the STTA method, were $r_\mu^L$ and  $r_\mu^R$ are the rank of the left and right DRMs $X_\mu$ and $Y_\mu$; see~\eqref{eq:def_Xmu_Ymu}. Since STTA and TT-HMT involve random quantities, we perform 30 trials and show the spread between the 20th and 80th percentile of the relative error $\|\mathcal T-\tilde{\mc T}\|_F / \|\tilde{\mc T}\|_F$ for the computed approximation $\mc T$.

A Python implementation of STTA and the code to reproduce  all the experiments is freely available at \url{https://github.com/RikVoorhaar/tt-sketch}.

\subsection{Baseline method TT-HMT}

We recall how the HMT method can be used for sketching tensors leading to what we call the TT-HMT method, summarized in \Cref{alg:hmt-sketch}. This algorithm produces the same approximation as the algorithms proposed in \cite{huberRandomizedTensorTrain2017, wolfLowRankTensor2019, cheRandomizedAlgorithmsApproximations2019}. An important numerical and implementation difference is that \Cref{alg:hmt-sketch} continues to operate on the original tensor $\mathcal{T}$ whereas the other algorithms successively reduce the size of $\mathcal{T}$, like the TT-SVD method. The benefit is that our formulation can  exploit structure in  $\mathcal{T}$ since only the original $\mathcal{T}$ needs to be accessed during the algorithm.

Superficially TT-HMT method is very similar to the STTA method. However, instead of performing a two-sided sketch using two independent DRMs $X_\mu$ and $Y_\mu$, the sketches $\Psi_\mu$ on line 4, 6 and 8 in \Cref{alg:hmt-sketch} are computed using only a right DRM $X_\mu$ and on the left we multiply by the partial contraction $C_{\leq \mu-1}$ of the TT-cores previously computed. Effectively, this partial contraction $C_{\leq \mu-1}$ can be viewed as a non-random left TT-DRM, and therefore the sketches $\Psi_\mu$ can be computed for any of the structured tensors mentioned in \Cref{sec:structured-sketches}. Note that the TT-HMT method is not streamable due to the orthogonalization on line 11, and requires at least $d$ passes of the data. Furthermore this means it is not possible to parallelize the for loop starting on line 1. On the other hand, the HMT method uses fewer and smaller sketches, making it faster than STTA for smaller rank approximations (for larger rank, this is no longer always true; see Section~\ref{sec:timing}), and it also produces orthogonalized TTs. It has a theoretical error bound in expectation with a smaller error constant than STTA; see  \cite{huberRandomizedTensorTrain2017}.

\begin{algorithm}
  \caption{TT-HMT}\label{alg:hmt-sketch}
  \begin{algorithmic}[1]
    \REQUIRE Tensor $\mathcal T$ of shape $n_1\times \dots \times n_d$; Right DRMs $X_\mu$ of size $(n_{\mu+1}\cdots n_d)\times r_\mu$.
    \ENSURE Cores $C_\mu$ of size $r_{\mu-1}\times n_\mu\times r_\mu$ for $1\leq \mu\leq d$ of a tensor train approximation of $\mc T$.
    
    \FOR{$\mu=1$ to $d$}
    \STATE $\mathcal T^{\leq \mu} \leftarrow \mathtt{reshape}(\mathcal T,\, (n_1\cdots n_\mu) \times (n_{\mu+1}\cdots n_d))$ 
    \IF{$\mu=1$}
      \STATE $\Psi_1 \leftarrow \mathcal T^{\leq 1}X_1$
      \ELSIF{$\mu<d$}
    \STATE $\Psi_\mu \leftarrow (C_{\leq \mu-1}^\top\tensor I_{n_\mu}) \mathcal T^{\leq \mu} X_\mu$
    \ELSE
    \STATE $\Psi_d \leftarrow (C_{\leq d-1}^\top\tensor I_{n_d}) \mathcal T^{\leq d}$
    \ENDIF
    \IF{$\mu<d$}
    \STATE $C_\mu^L \leftarrow \mathtt{orth}(\Psi_\mu^L)$ \hfill // Q factor of QR decomposition
    \STATE $C_{\leq \mu}\leftarrow (C_{\leq \mu-1}\otimes I_{n_{\mu}})C_\mu^L$ \hfill // Not formed explicitly
    \STATE $C_\mu \leftarrow \mathtt{reshape}(C_\mu^L,\,r_{\mu-1}^L\times n_\mu \times r_\mu^L)$
    \ELSE
    \STATE $C_d\leftarrow \Psi_d$
    \ENDIF
    \ENDFOR
  \end{algorithmic}
\end{algorithm}

\subsection{Sketching dense tensors}

\begin{figure}[htb]
  \centering
  \includegraphics[width=0.7\textwidth]{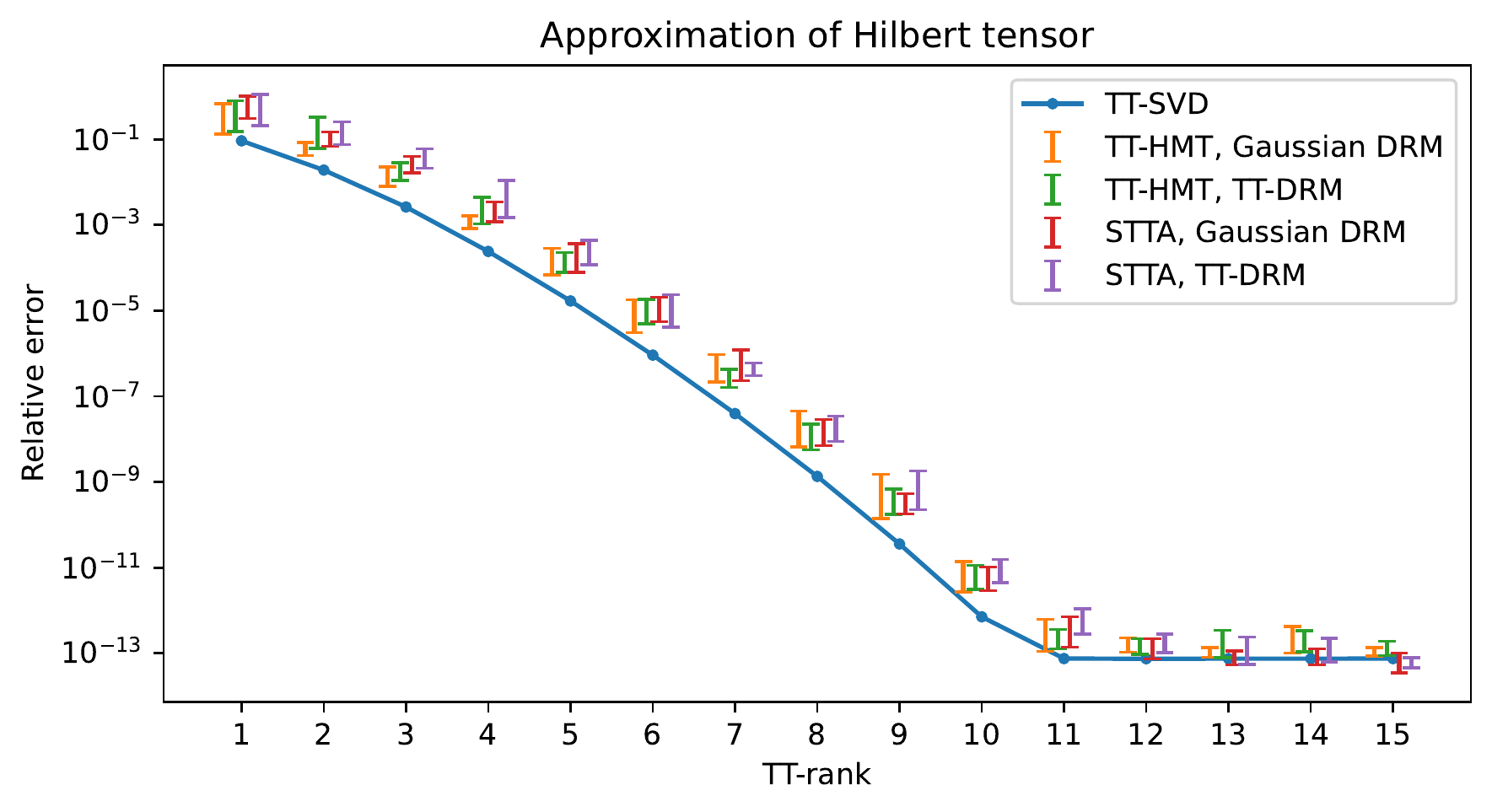}
  \caption{Approximation errors for the Hilbert tensor different TT-ranks for the TT-SVD, STTA and TT-HMT methods. The markers are ordered left-to-right in the same order as the legend.}\label{fig:plot-hilbert}
\end{figure}

\begin{figure}[htb]
  \centering
  \includegraphics[width=0.7\textwidth]{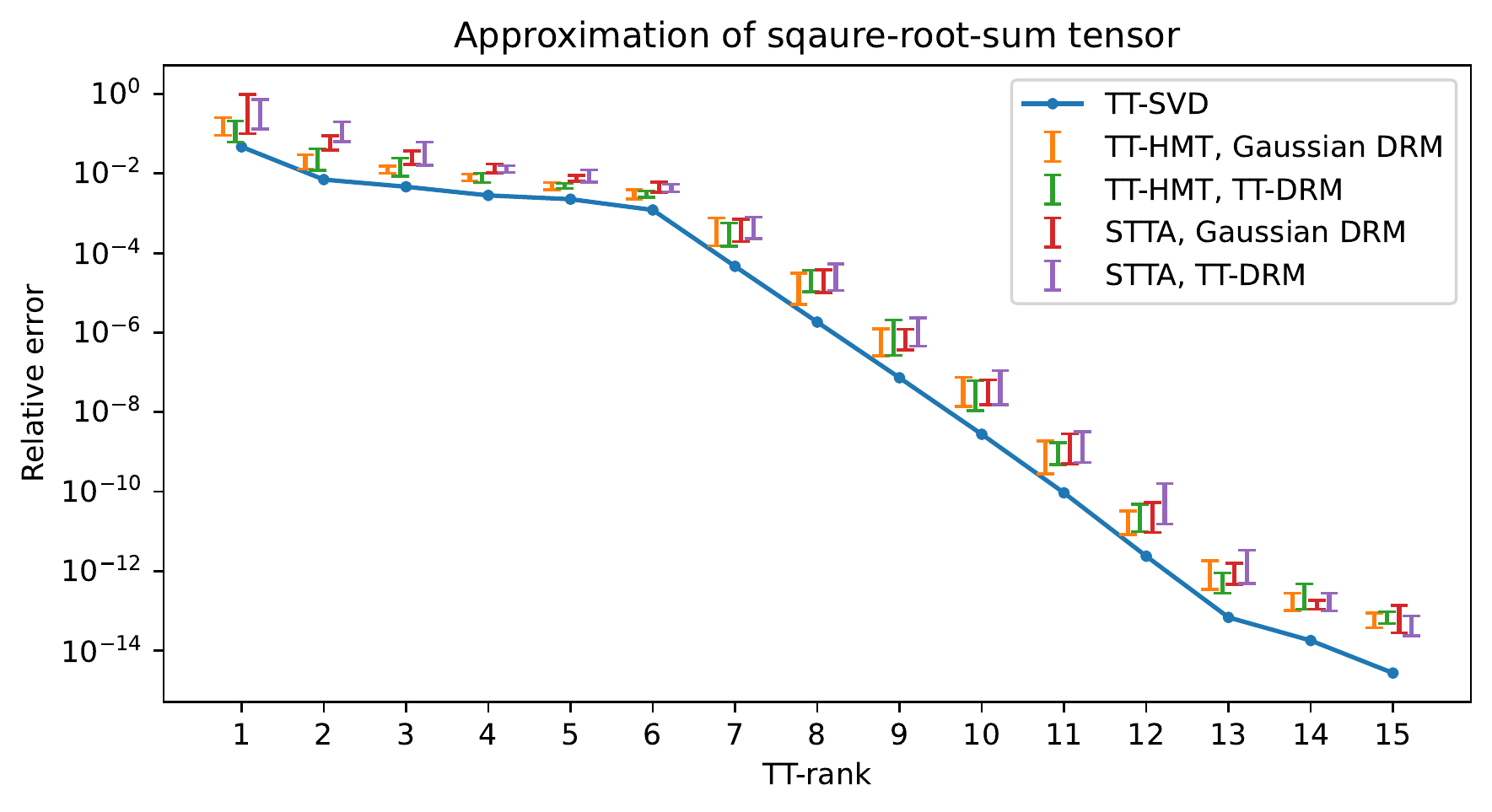}
  \caption{Same setting as in \cref{fig:plot-hilbert} but for the square-root-sum tensor.}\label{fig:plot-sum-sqrt}
\end{figure}

We first consider the problem of sketching dense tensors of shape $n_1 \times \cdots \times n_d$. In particular, we will test the performance of the algorithms on the Hilbert tensor 
\begin{equation}
  \mathcal T_{\mathrm{Hilbert}}[i_1,\dots,i_d] = (i_1+\cdots+i_d - d+1)^{-1}
\end{equation}
and the square-root-sum tensor
\begin{equation}
  \mathcal T_{\mathrm{sqrt}}[i_1,\dots,i_d] = \sqrt{\sum_{j=1}^d \frac{n_j-i_j}{n_j-1}a + \frac{i_j-1}{n_j-1}b},
\end{equation}
with $0<a<b$. 
Both tensors have exponential decay of the singular values of the unfoldings $\mc T^{\leq \mu}$ for all $\mu$; see, e.g.,~\cite{hackbuschTensorSpacesNumerical2019}. For $\mathcal T_{\mathrm{sqrt}}$, the decay worsens as $a$ approaches zero. In our experiments we take $a=0.2$ and $b=2$. Furthermore, we take $d=7$ and $n_1 = \cdots = n_7 = 5$ for $\mathcal T_{\mathrm{Hilbert}}$, and $d=5$ and $n_1=\dots=n_5=10$ for $\mathcal T_{\mathrm{sqrt}}$. The approximation errors for different TT ranks are shown in \cref{fig:plot-hilbert} and \cref{fig:plot-sum-sqrt}

We observe that the error of STTA and TT-HMT are comparable and only slightly worse than TT-SVD, and has acceptable variance. For both STTA and TT-HMT, there is not much difference between the use of Gaussian and TT DRMs when comparing approximation error and its variance. This is surprising since TT DRMs have less randomness than Gaussian DRMs.

\subsection{Sketching a sum of TT and a sparse tensor}
One of the major advantages of STTA is that it can approximate the sum of several tensors in parallel, possibly in different formats. As an example we consider the problem of approximating the sum of a low-rank TT and a sparse tensor. 

\begin{figure}[htb]
  \centering
  \includegraphics[width=0.8\textwidth]{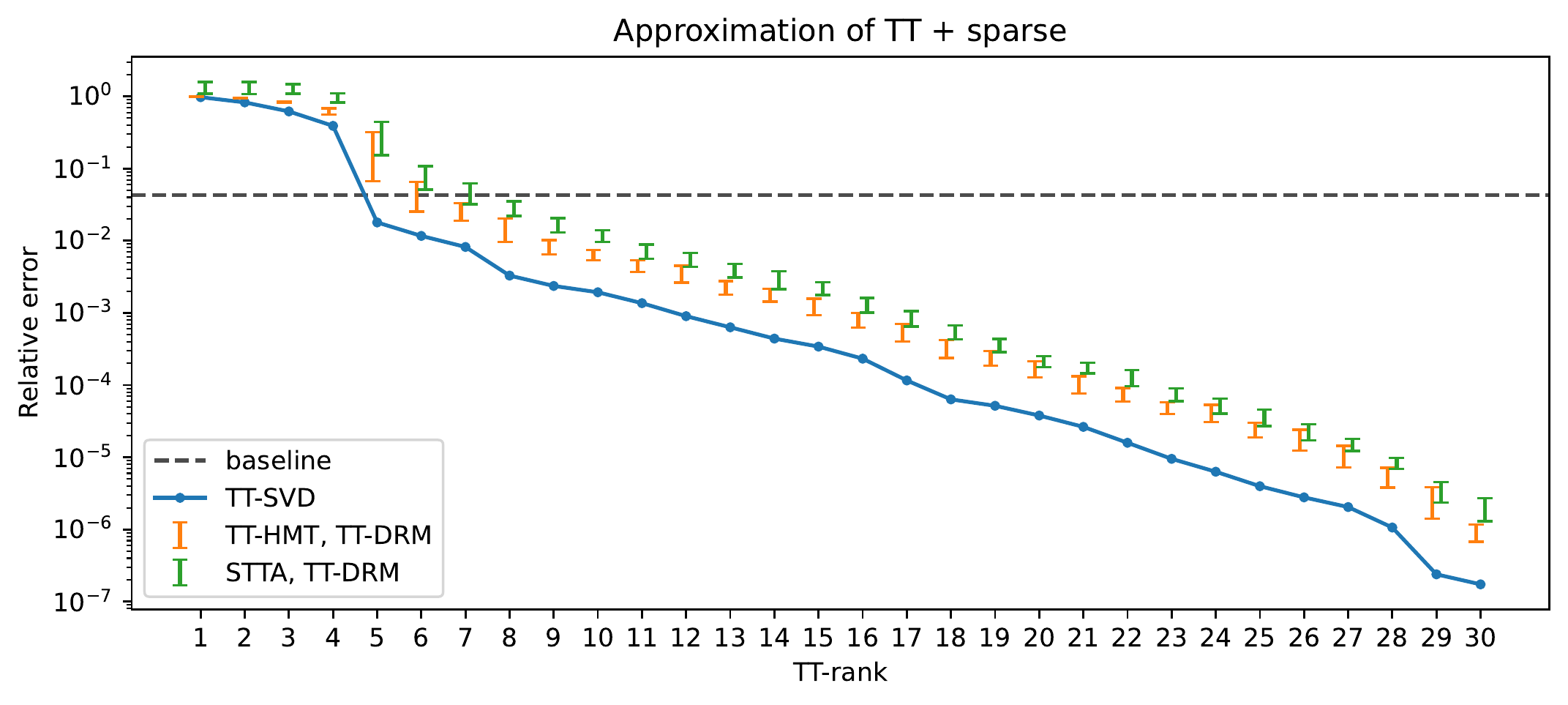}
  \caption{Approximation errors for a sum of a low-rank TT and a sparse tensor for different TT-ranks for the TT-SVD, STTA (right markers) and TT-HMT (left markers) methods. As baseline we show the error obtained when ignoring the sparse tensor and just using the low-rank TT as approximation.}\label{fig:plot-tt-plus-sparse}
\end{figure}

We consider a rank-5 TT of size $(10,10,10,10,10)$ consisting of cores with i.i.d. $\mathrm N(0,1/r^2)$ entries. To this we add a random sparse tensor with 100 nonzero entries. The nonzero entries are normally distributed with standard deviation ranging log-uniformly between $10^{-3}$ and $10^{-20}$. The resulting approximation is shown in \cref{fig:plot-tt-plus-sparse}.

\subsection{Sketching a sum of TT}\label{sec:experiment-sum-of-tt}
\begin{figure}[htb]
  \centering
  \includegraphics[width=0.8\textwidth]{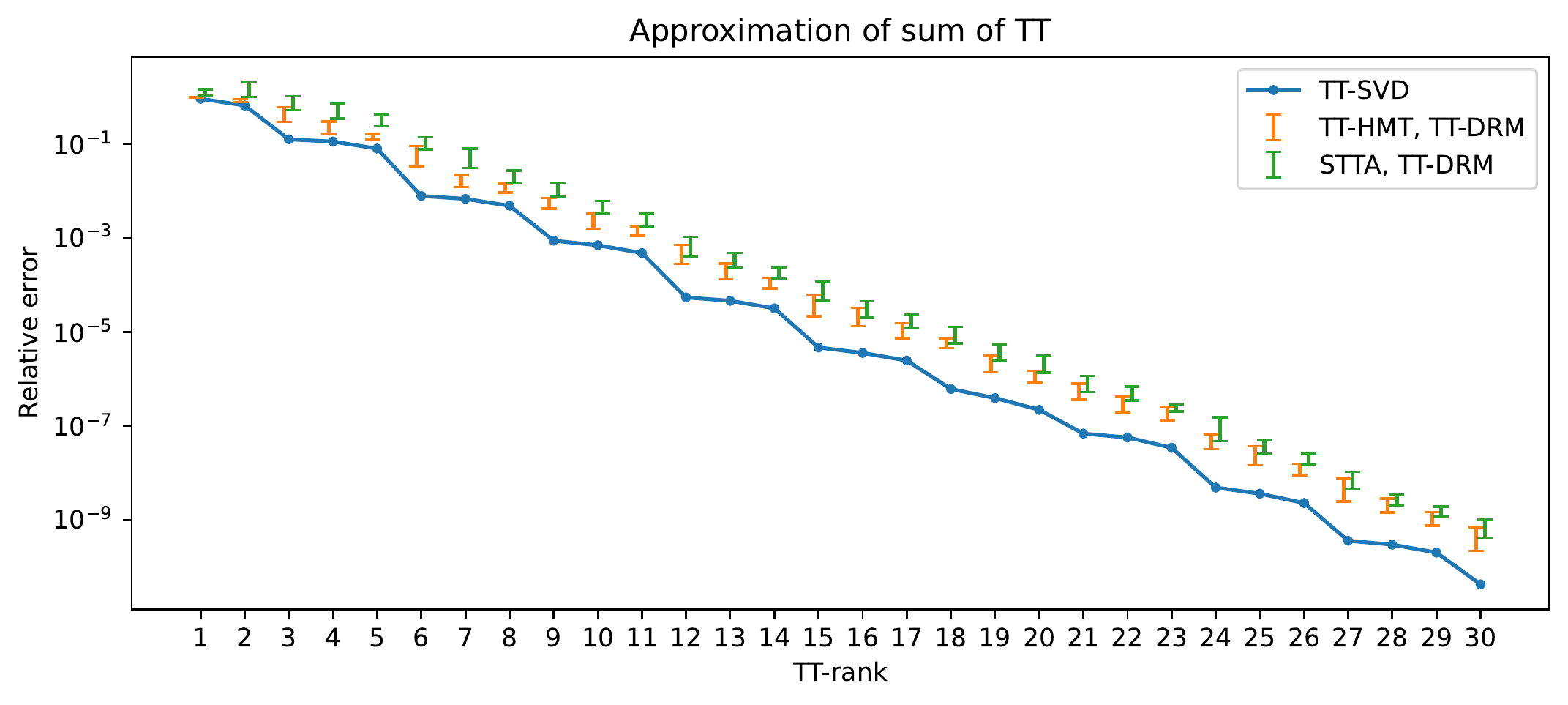}
  \caption{Approximation errors for a sum of a low-rank TTs for different TT-ranks for the TT-SVD, STTA (right) and TT-HMT (left) algorithms. }\label{fig:plot-tt-sum}
\end{figure}

Next we consider the problem of sketching a sum of low-rank TTs. Specifically we take a sum of 20 rank-3 TTs of size $10\times 10\times 10\times 10\times 10$ of form
\[
\mc T = \sum_{i=0}^{19} 10^{-i} \mc T_i,
\]
where each $\mc T_i$ is a random TT of rank 3 consisting of cores with i.i.d. $\mathrm N(0,1/r^2n)$ entries. The resulting approximation errors for the different algorithms are shown below in \cref{fig:plot-tt-sum}.

\subsection{Sketching a CP tensor}
\begin{figure}[htb]
  \centering
  \includegraphics[width=0.8\textwidth]{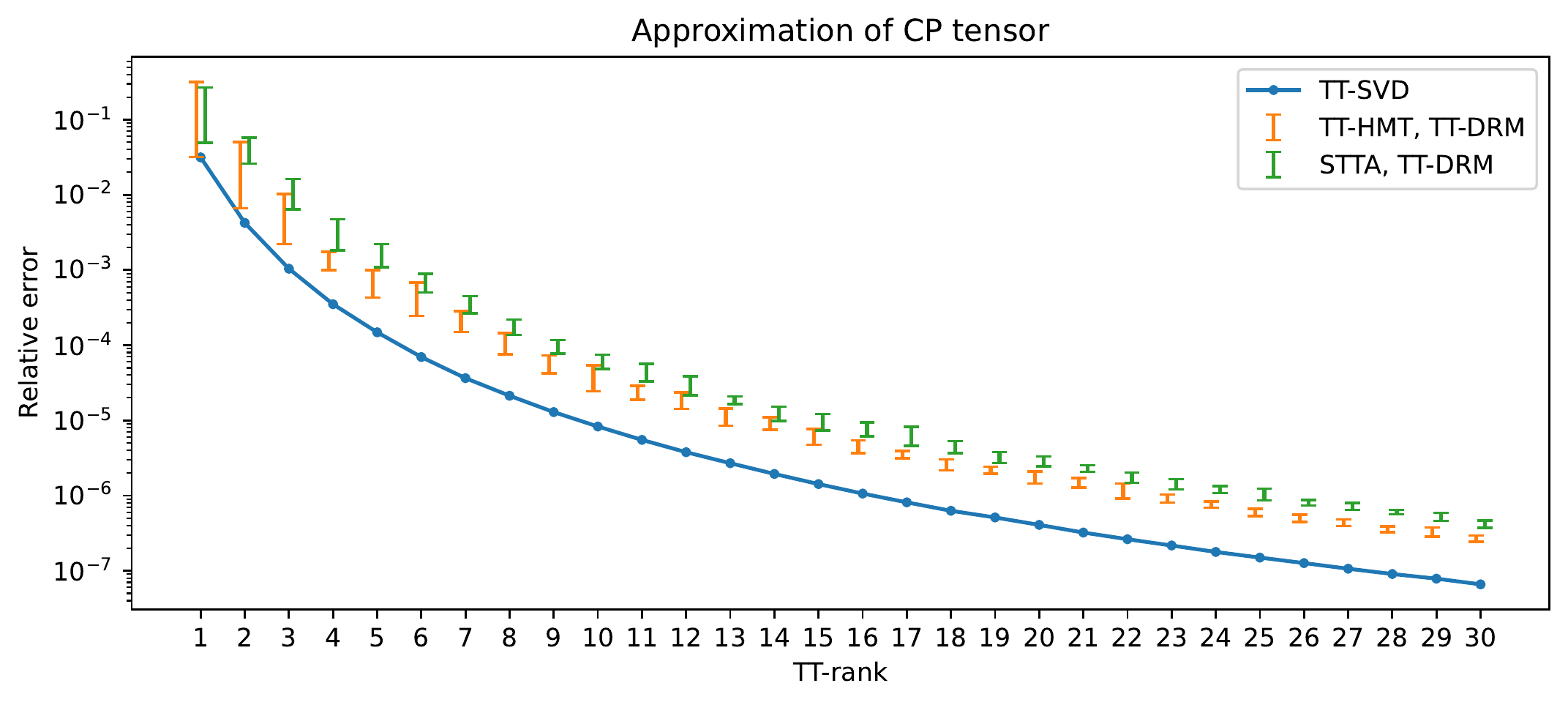}
  \caption{Same setting as in \cref{fig:plot-tt-sum} but for a CP tensor. }\label{fig:plot-cp}
\end{figure}

We consider the problem of sketching a CP tensor. The CP tensor is of form
\[
  \mathcal T = \sum_{i=1}^{100} \frac{1}{i^5} v_{1,i}\tensor\cdots \tensor v_{5,i},
\]
where each $v_{\mu,i}\in \R^{10}$ is a normalized Gaussian vector.  The results are shown in \cref{fig:plot-cp}.

\subsection{The effect of oversampling}
\begin{figure}[htb]
  \centering
  \includegraphics[width=0.7\textwidth]{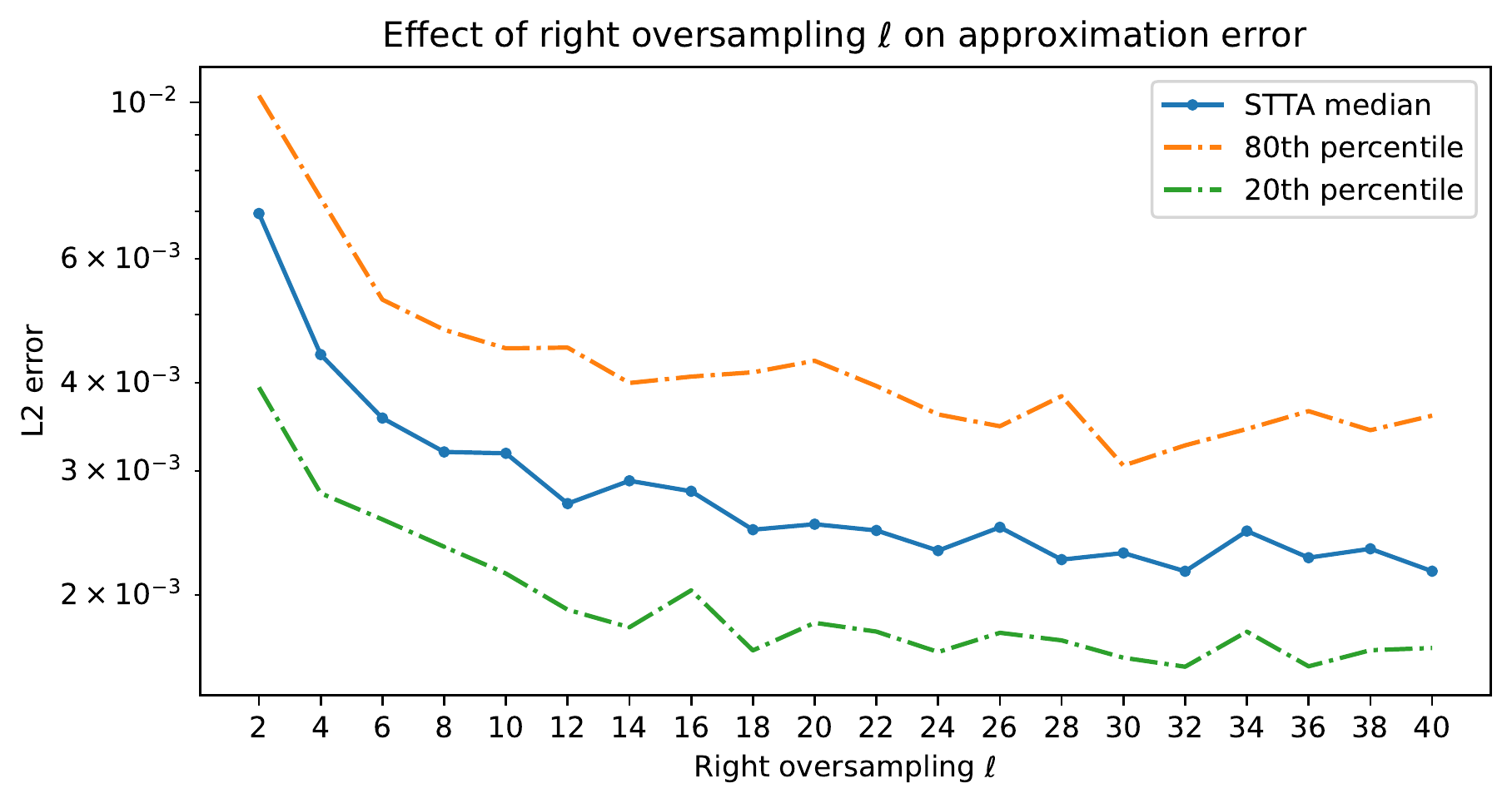}
  \caption{Effect of oversampling parameter $\ell_\mu=|r_\mu^L-r_\mu^R|$ on the approximation error for STTA.}\label{fig:plot-right-oversampling}
\end{figure}

The rank $r_\mu^L$ of the left DRMs $X_\mu$ is required to be smaller than the rank $r_\mu^R$ of the left DRMs $Y_\mu$ for all $\mu$ (or vice versa). Let $r_\mu^R=r_\mu^L+\ell$. We investigate the effect of varying the oversampling $\ell>0$ on the approximation error. As test problem we will use the same sum of TTs as in \cref{sec:experiment-sum-of-tt}, and we fix the TT rank of the approximation at $r_\mu = r=10$ for all $\mu$. The result of this experiment is shown below in \cref{fig:plot-right-oversampling}

Observe that increasing the oversampling $\ell$ improves both the median and the variance of the approximation error. (Note that the scale for the error in \cref{fig:plot-right-oversampling} is logarithmic.)  However, in this particular case there are quickly diminishing returns after $\ell \approx 10$. Since the cost of all algorithms depends at most linearly on $\ell$,  setting $\ell=r$ (as we have done in the experiments above), will only double the computational cost compared to no oversampling at all. This provides a reasonable trade-off between approximation error and speed.

\subsection{Dependence of the error on the tensor order}\label{sec:experiment-tensor-order}

\Cref{thm:error-bound} implies that the approximation constant of the STTA method with Gaussian DRMs scales exponentially with the number of modes $d$. In practice, it appears that the approximation constant scales as $\sqrt{d-1}$ for Gaussian and TT DRMs, like with the classical TT-SVD and the TT-HMT method using Gaussian DRMs. To test this, we compute a rank 10 approximation of a TT of uniform TT-rank $r=30$ and size $30^{d}$ for different values of $d$. For each unfolding, the TT has singular values decaying exponentially between $\sqrt{r}$ and $\sqrt{r}\cdot 10^{-20}$ (this keeps the Frobenius norm of the TT approximately constant as $d\to\infty$). 

\begin{figure}[htb]
    \centering
    \includegraphics[width=0.65\textwidth]{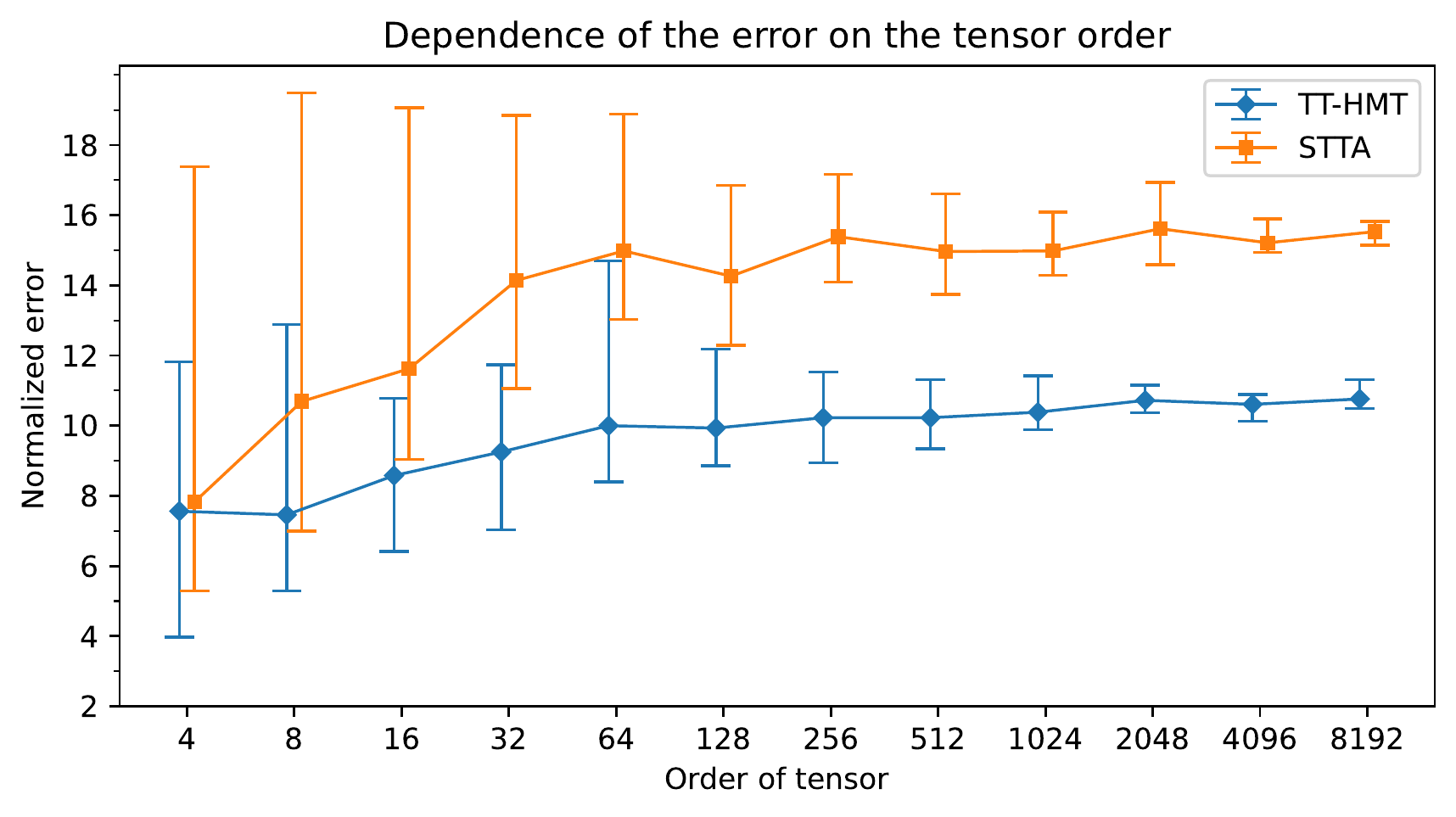}
    \caption{Normalized approximation errors for TTs with varying tensor orders for rank 10 approximations obtained with TT-HMT (left markers) and STTA (right markers). The approximation error is normalized by dividing by the approximation error of the TT-SVD approximation for the same rank. }
    \label{fig:plot-dimension-scaling}
\end{figure}

In \cref{fig:plot-dimension-scaling} we compare the scaling of the approximation error with the tensor order $d$ for TT-HMT and STTA using TT DRMs. We normalize the approximation error by dividing by the approximation error of the TT-SVD method for the same rank. Note that in this context by TT-SVD we mean the standard SVD-based rounding procedure for TTs.

We observe that the STTA method converges to an approximation error roughly 13 times worse than the TT-SVD method, whereas the HMT method converges to an approximation error roughly 8 times worse than the TT-SVD method. This suggests that the approximation constant of the STTA method asymptotically scales as $\sqrt{d-1}$, just like the TT-SVD and TT-HMT methods.

We remark that the scaling $1/r_{\mu}^L$ mentioned in \Cref{def:tt-drm} of TT-DRMs becomes essential for the large number of modes used in this experiment. Different scalings of the cores of the TT-DRMs leads to severe problems with underflow or overflow.

\subsection{Timing comparison}\label{sec:timing}
\begin{figure}[htb]
\begin{minipage}{.49\textwidth}
  \centering
  \includegraphics[width=1\textwidth]{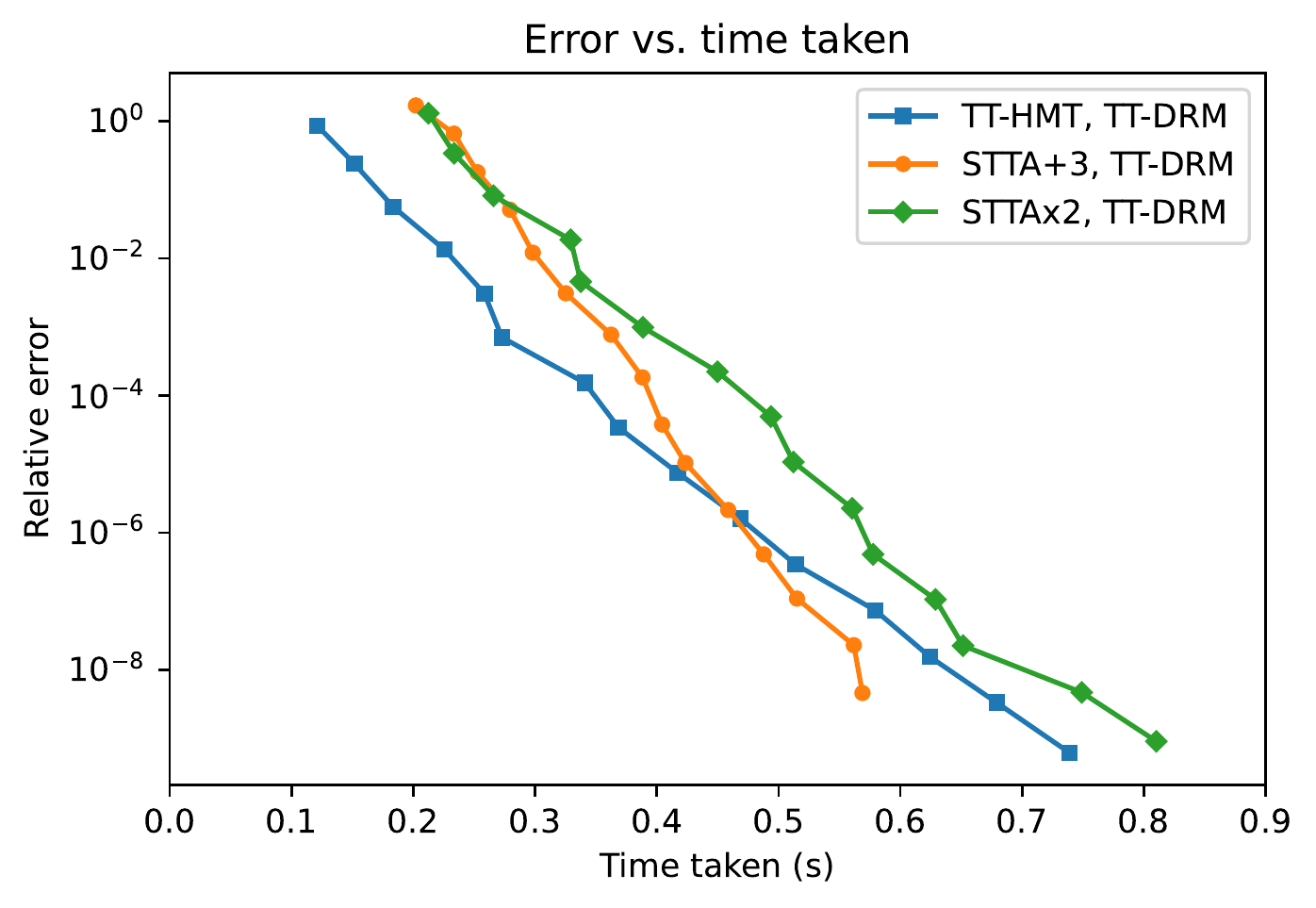}
\end{minipage}
\hfill
\begin{minipage}{.49\textwidth}
  \centering
  \includegraphics[width=1\textwidth]{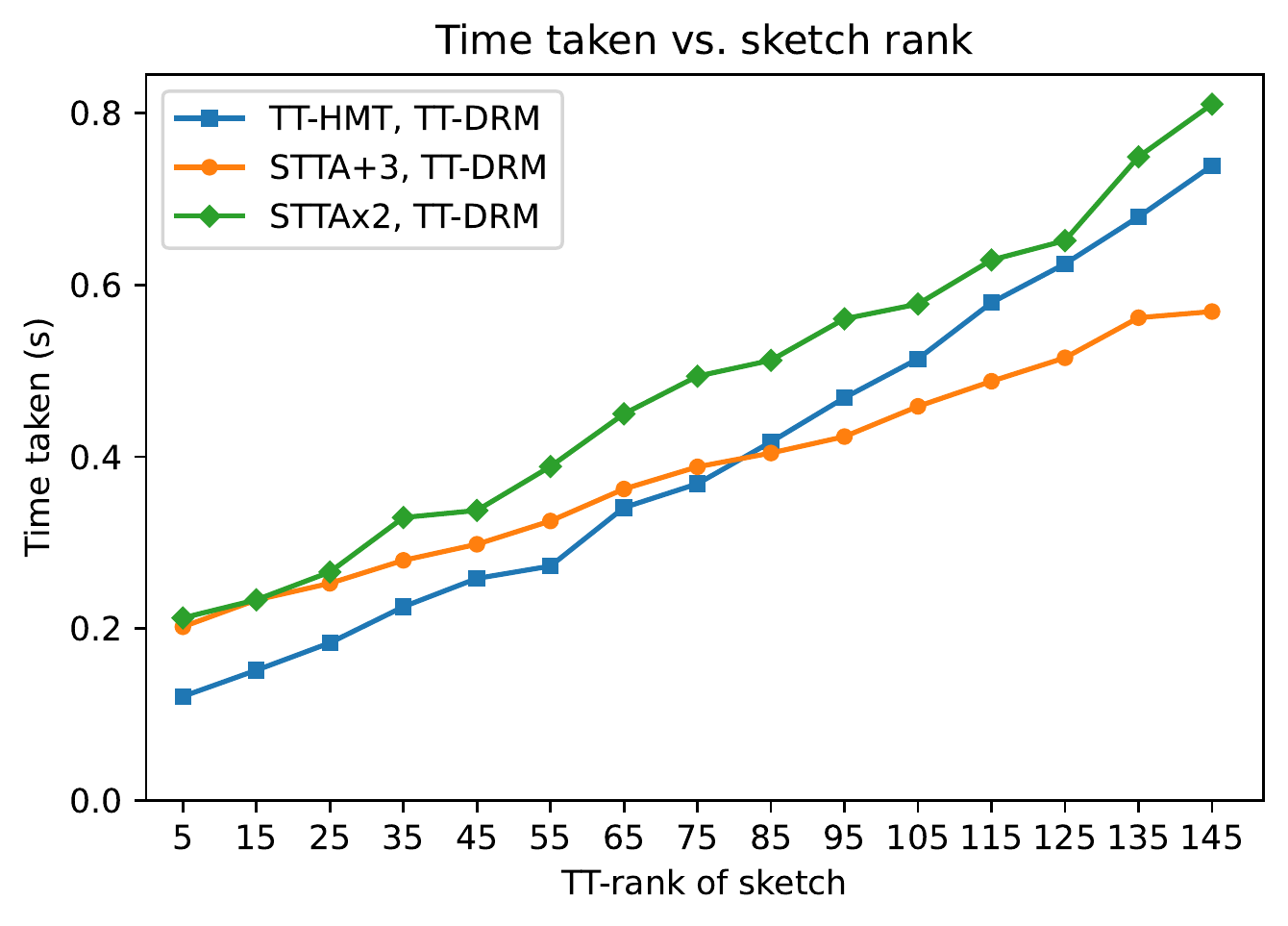}
\end{minipage}
  \caption{Relation between timing and error (left) and the TT-rank of the TT-DRM and time taken (right) for TT-HMT, and STTA with respectively $r_\mu^L=r_\mu^R+3$ (STTA+3) and $r_\mu^L=2r_\mu^R$ (STTAx2).}\label{fig:plot-timings}
\end{figure}

We consider the trade-off between speed, rank and approximation error for STTA and the TT-HMT method. In particular we approximate a random TT of size $150\times 150\times 150\times 150\times 150$ of TT-rank $150$. The singular values of each unfolding were made to decay exponentially between 1 and and $10^{-10}$. In \cref{fig:plot-timings} we report the time taken and relative error of approximation for the TT-HMT and STTA method using TT-DRMs of varying ranks. All experiments were performed on a 2x64 core AMD Epyc 7742. For the approximation error and the time taken we report the median over 100 trials.

While the TT-HMT method is faster for a low-rank approximation, the STTA method is faster than TT-HMT for sufficiently large rank approximations. This is because the cost of the orthogonalization step in TT-HMT becomes dominant for sufficiently large ranks, and offsets the cost of the second sketch in STTA. This is especially true if we use a low amount of oversampling ($r_\mu^L=r_\mu^R+3$).

\section{Conclusions}

This paper introduced Streaming Tensor Train Approximation (STTA) -- an algorithm for approximating a tensor in the tensor train format that is both streamable and one-pass. Since the method computes the approximation using only two-sided sketches of the given tensor, STTA can efficiently update the approximation after linear updates without needing to access the original tensor. In addition, STTA can be applied to specific structures of the tensor such as (linear combinations of) sparsity and various low-rank tensor formats (CP, Tucker, tensor train). In case of sketching with Gaussian dimension reduction matrices and with sufficiently large oversampling, we proved a quasi-optimal error bound in expectation by extending existing results from the Generalized Nyström method for matrices. In addition, the numerical experiments showed that the approximation error concentrates relatively quickly around its mean without needing much oversampling in practice. This is also true when sketching using structured dimension reduction matrices based on tensor trains. Due to the high-level of parallelism, STTA is also faster than computing tensor train approximations based on the Halko--Martinson--Tropp method in a multi-core environment.

\bibliography{main}
\bibliographystyle{hapalike}

\end{document}